\input amstex
\documentstyle{amsppt}
%
\catcode`@=11
\redefine\output@{%
  \def\break{\penalty-\@M}\let\par\endgraf
  \ifodd\pageno\global\hoffset=105pt\else\global\hoffset=8pt\fi  
  \shipout\vbox{%
    \ifplain@
      \let\makeheadline\relax \let\makefootline\relax
    \else
      \iffirstpage@ \global\firstpage@false
        \let\rightheadline\frheadline
        \let\leftheadline\flheadline
      \else
        \ifrunheads@ 
        \else \let\makeheadline\relax
        \fi
      \fi
    \fi
    \makeheadline \pagebody \makefootline}%
  \advancepageno \ifnum\outputpenalty>-\@MM\else\dosupereject\fi
}
\catcode`\@=\active
\nopagenumbers
\def\negskp{\hskip -2pt}
\def\Alpha{\operatorname{A}}
\def\MatGrSO{\operatorname{SO}}
\def\MatGrSL{\operatorname{SL}}
\def\compos{\,\raise 1pt\hbox{$\sssize\circ$} \,}
\accentedsymbol\bd{\kern 2pt\bar{\kern -2pt d}}
\accentedsymbol\bbd{\kern 2pt\bar{\kern -2pt\bold d}}
\accentedsymbol\bundleboldX{\raise 8pt \vbox{\hsize=7pt\noindent
  \vrule height 0.4pt depth 3pt width 1pt\vbox{\hrule width 5pt}}
  \kern -5pt\bold X}
\def\msum{\operatornamewithlimits{\sum^3\!{\ssize\ldots}\!\sum^3}}
\def\vtrule{\vrule height 12pt depth 6pt}
\def\vtttrule{\vrule height 12pt depth 19pt}
\def\boxit#1#2{\vcenter{\hsize=122pt\offinterlineskip\hrule
  \line{\vtttrule\hss\vtop{\hsize=120pt\centerline{#1}\vskip 5pt
  \centerline{#2}}\hss\vtttrule}\hrule}}
\def\blue#1{#1}
\def\red#1{#1}
\catcode`#=11\def\diez{#}\catcode`#=6
\catcode`_=11\def\podcherkivanie{_}\catcode`_=8
\def\mycite#1{\cite{\blue{#1}}\immediate\special{ps:
     ShrHPSdict begin /ShrBORDERthickness 0 def}}
\def\myciterange#1#2#3#4{\cite{\blue{#2#3#4}}\immediate\special{ps:
     ShrHPSdict begin /ShrBORDERthickness 0 def}}
\def\mytag#1{%
    \tag#1}
\def\mythetag#1{\thetag{\blue{#1}}\immediate\special{ps:
     ShrHPSdict begin /ShrBORDERthickness 0 def}}
\def\myrefno#1{\no#1}
\def\myhref#1#2{\blue{#2}\immediate\special{ps:
     ShrHPSdict begin /ShrBORDERthickness 0 def}}
\def\myEarXivlink{\myhref{http://arXiv.org}{http:/\negskp/arXiv.org}}
\def\myGeoCities{\myhref{http://www.geocities.com}{GeoCities}}
\def\mytheorem#1{\csname proclaim\endcsname{Theorem #1}}
\def\mythetheorem#1{\blue{#1}\immediate\special{ps:
     ShrHPSdict begin /ShrBORDERthickness 0 def}}
\def\mylemma#1{\csname proclaim\endcsname{Lemma #1}}

\def\mycorollary#1{\csname proclaim\endcsname{Corollary #1}}

\def\mydefinition#1{\definition{Definition #1}}
\def\mythedefinition#1{\blue{#1}\immediate\special{ps:
     ShrHPSdict begin /ShrBORDERthickness 0 def}}

\pagewidth{360pt}
\pageheight{606pt}
\topmatter
\title
A note on Kosmann-Lie derivatives of Weyl spinors.
\endtitle
\author
R.~A.~Sharipov
\endauthor
\address 5 Rabochaya street, 450003 Ufa, Russia\newline
\vphantom{a}\kern 12pt Cell Phone: +7(917)476 93 48
\endaddress
\email \vtop to 30pt{\hsize=280pt\noindent
\myhref{mailto:r-sharipov\@mail.ru}
{r-sharipov\@mail.ru}\newline
\myhref{mailto:R\podcherkivanie Sharipov\@ic.bashedu.ru}
{R\_\hskip 1pt Sharipov\@ic.bashedu.ru}\vss}
\endemail
\urladdr
\vtop to 20pt{\hsize=280pt\noindent
\myhref{http://www.geocities.com/r-sharipov}
{http:/\negskp/www.geocities.com/r-sharipov}\newline
\myhref{http://www.freetextbooks.boom.ru/index.html}
{http:/\negskp/www.freetextbooks.boom.ru/index.html}\vss}
\endurladdr
\abstract
    Kosmann-Lie derivatives in the bundle of Weyl spinors are considered.
It is shown that the basic spin-tensorial fields of this bundle are 
constants with respect to these derivatives.
\endabstract
\subjclassyear{2000}
\subjclass 53B30, 81T20, 22E70\endsubjclass
\endtopmatter
\loadbold
\loadeufb
\TagsOnRight
\document
\accentedsymbol\hatboldsymbolgamma{\kern 1.3pt\hat{\kern -1.3pt
   \boldsymbol\gamma}}
\accentedsymbol\hatgamma{\kern 1.3pt\hat{\kern -1.3pt\gamma}}

\head
1. Introduction. 
\endhead
    Lie derivatives arise in studying continuous symmetries of various 
geometric structures on manifolds. They are also used in symmetry analysis 
of ordinary and partial differential equations (see \mycite{1}). In 
general relativity the bundle of Weyl spinors $SM$ is a special geometric 
structure built over the space-time manifold $M$. The main goal of this 
paper is to clarify the procedure of applying Lie derivatives to the basic
attributes of this geometric structure, i\.\,e\. to the basic spin-tensorial 
fields associated with the bundle of Weyl spinors.
\head
2. Lie derivatives of spatial structures.
\endhead
     Let $M$ be a space-time manifold of general relativity.
This means that it is a four-dimensional orientable manifold 
equipped with a Minkowski type metric $\bold g$ and with a 
polarization. A polarization, which is typically not mentioned,
is a geometric structure that marks the future half light cone 
in the tangent space $T_p(M)$ for each point $p\in M$ (see more 
details in \mycite{2}). A Lie derivative $L_{\bold X}$ is
usually given by some vector field $\bold X$ in $M$. Once such
a vector field $\bold X$ is fixed, it produces a one-parametric 
local group of local diffeomorphisms:
$$
\hskip -2em
\varphi_{\varepsilon}\!:\,M\to M.
\mytag{2.1}
$$
The letter $t$ is typically used for the parameter of this local group (see 
\mycite{3}), but here we use the Greek letter $\varepsilon$ since $t$ in 
physics is reserved for the time variable. The local diffeomorphisms 
\mythetag{2.1} induce the local diffeomorphisms 
$$
\xalignat 2
\hskip -2em
\varphi_{\varepsilon*}\!:\,TM\to TM,
&&\varphi^*_{-\varepsilon}\!:\,T^*\!M\to T^*\!M
\mytag{2.2}
\endxalignat 
$$
in tangent and cotangent bundles respectively. These induced 
diffeomorphisms \mythetag{2.2} act as linear mappings in fibers 
of $TM$ and $T^*\!M$. For this reason they can be extended to 
local diffeomorphisms of tensor bundles:
$$
\hskip -2em
\varphi_{\varepsilon}\!:\,T^r_sM\to T^r_sM.
\mytag{2.3}
$$
Here in \mythetag{2.3} through $T^r_sM$ we denote the following tensor
product of $r$ copies of the tangent bundle $TM$ and $s$ copies of the
cotangent bundle $T^*\!M$:
$$
\hskip -2em
T^r_sM=\overbrace{TM\otimes\ldots\otimes TM}^{\text{$r$ times}}
\otimes\underbrace{T^*\!M\otimes\ldots\otimes T^*\!M}_{\text{$s$ times}}.
\mytag{2.4}
$$\par
     Let's study the diffeomorphisms \mythetag{2.1} and \mythetag{2.2} 
in more details. Assume that $p$ and $q$ are two points of the space-time
manifold $M$ such that $q=\varphi_{\varepsilon}(p)$. Then 
$p=\varphi_{-\varepsilon}(q)$ and we have the following commutative 
diagram:
$$
\hskip -2em
\CD
T_p(M)@>\dsize\varphi_{\varepsilon*}>> T_q(M)\\
@V\pi VV @VV\pi V\\
p@>\dsize\varphi_{\varepsilon}>>q\\
@A\pi AA @AA\pi A\\
T^*_p(M)@>\dsize\varphi^*_{-\varepsilon}>> T^*_q(M).
\endCD
\mytag{2.5}
$$
Assume that we have some local chart with the coordinates $x^0,\,x^1,\,x^2,
\,x^3$ in some neighborhood of the point $q$. Assume also that $\varepsilon$
is small enough so that the point $p=\varphi_{-\varepsilon}(q)$ in 
\mythetag{2.5} is covered by the same local chart. Then the coordinates 
of the points $p$ and $q$ in this chart are related to each other as follows:
$$
\left\{\aligned
&x^0=u^0(\varepsilon,y^0,y^1,y^2,y^3),\\
&x^1=u^1(\varepsilon,y^0,y^1,y^2,y^3),\\
&x^2=u^2(\varepsilon,y^0,y^1,y^2,y^3),\\
&x^3=u^3(\varepsilon,y^0,y^1,y^2,y^3),
\endaligned\right.
\qquad\quad
\left\{\aligned
&y^0=u^0(-\varepsilon,x^0,x^1,x^2,x^3),\\
&y^1=u^1(-\varepsilon,x^0,x^1,x^2,x^3),\\
&y^2=u^2(-\varepsilon,x^0,x^1,x^2,x^3),\\
&y^3=u^3(-\varepsilon,x^0,x^1,x^2,x^3).
\endaligned\right.\quad
\mytag{2.6}
$$
Using \mythetag{2.6}, we define the matrices $\Phi(\varepsilon)$ and 
$\Phi(-\varepsilon)$ with the components
$$
\Phi^i_{\!j}(\varepsilon)=\frac{\partial u^i(\varepsilon,
y^0,y^1,y^2,y^3)}{\partial y^j},
\qquad
\Phi^i_{\!j}(-\varepsilon)=\frac{\partial u^i(-\varepsilon,
x^0,x^1,x^2,x^3)}{\partial x^j}.
\quad
\mytag{2.7}
$$\par
     Let $Y^{i_1\ldots\,i_r}_{j_1\ldots\,j_s}(x^0,x^1,x^2,x^3)$ 
be the components of some tensorial field $\bold Y$ of the type 
$(r,s)$ in the local coordinates $x^0,\,x^1,\,x^2,\,x^3$ and let
$\varphi_{\varepsilon}(\bold Y)^{i_1\ldots\,i_r}_{j_1\ldots
\,j_s}(x^0,x^1,x^2,x^3)$ be the components of its image 
$\varphi_{\varepsilon}(\bold Y)$ under the mapping \mythetag{2.3}. 
Then we have
$$
\gathered
\varphi_{\varepsilon}(\bold Y)^{i_1\ldots\,i_r}_{j_1\ldots
\,j_s}(x^0,x^1,x^2,x^3)=
\msum\Sb h_1,\,\ldots,\,h_r\\k_1,\,\ldots,\,k_s\endSb 
\Phi^{i_1}_{\!h_1}(\varepsilon)
\,\ldots\,\Phi^{i_r}_{\!h_r}(\varepsilon)\,\times\\
\vspace{1ex}
\qquad\qquad\qquad\qquad
\times\,\Phi^{k_1}_{\!j_1}(-\varepsilon)
\,\ldots\,\Phi^{k_s}_{\!j_s}(-\varepsilon)
\ Y^{h_1\ldots\,h_r}_{k_1\ldots\,k_s}(y^0,y^1,y^2,y^3).
\endgathered\qquad\qquad
\mytag{2.8}
$$
According to \mycite{3}, the Lie derivative $L_{\bold X}$ 
applied  to $\bold Y$ is defined as follows:
$$
\hskip -2em
L_{\bold X}(\bold Y)=\lim_{\varepsilon\to\,0}
\frac{\bold Y-\varphi_{\varepsilon}(\bold Y)}{\varepsilon}
=-\frac{d\varphi_{\varepsilon}(\bold Y)}{d\varepsilon}
\hbox{\vrule height14pt depth8pt width 0.5pt}_{\ \varepsilon=0}.
\mytag{2.9}
$$
Let $X^0,\,X^1,\,X^2,\,X^3$ be the components of the vector field 
$\bold X$ in the local coordinates $x^0,\,x^1,\,x^2,\,x^3$. Then
for $\varepsilon\to 0$ we have the following Taylor expansions of 
the functions \mythetag{2.6} representing the mappings
$\varphi_{\varepsilon}$ and $\varphi_{-\varepsilon}$:
$$
\align
&\hskip -2em
\left\{\aligned
&u^0(\varepsilon,y^0,y^1,y^2,y^3)=y^0+X^0(y^0,y^1,y^2,y^3)
\,\varepsilon+\,\ldots,\\\
&u^1(\varepsilon,y^0,y^1,y^2,y^3)=y^1+X^1(y^0,y^1,y^2,y^3)
\,\varepsilon+\,\ldots,\\
&u^2(\varepsilon,y^0,y^1,y^2,y^3)=y^2+X^2(y^0,y^1,y^2,y^3)
\,\varepsilon+\,\ldots,\\
&u^3(\varepsilon,y^0,y^1,y^2,y^3)=y^3+X^3(y^0,y^1,y^2,y^3)
\,\varepsilon+\,\ldots,
\endaligned\right.
\mytag{2.10}\\
\vspace{2ex}
&\hskip -2em
\left\{\aligned
&u^0(-\varepsilon,x^0,x^1,x^2,x^3)=x^0-X^0(x^0,x^1,x^2,x^3)
\,\varepsilon+\,\ldots,\\\
&u^1(-\varepsilon,x^0,x^1,x^2,x^3)=x^1-X^1(x^0,x^1,x^2,x^3)
\,\varepsilon+\,\ldots,\\
&u^2(-\varepsilon,x^0,x^1,x^2,x^3)=x^2-X^2(x^0,x^1,x^2,x^3)
\,\varepsilon+\,\ldots,\\
&u^3(-\varepsilon,x^0,x^1,x^2,x^3)=x^3-X^3(x^0,x^1,x^2,x^3)
\,\varepsilon+\,\ldots.
\endaligned\right.
\mytag{2.11}
\endalign
$$
Applying \mythetag{2.10} and \mythetag{2.11} to \mythetag{2.7}, we 
derive
$$
\align
&\hskip -2em
\Phi^i_{\!j}(\varepsilon)=\delta^i_j+\frac{\partial X^i(x^0,x^1,x^2,x^3)}
{\partial x^j}\ \varepsilon+\,\ldots,
\mytag{2.12}\\
&\hskip -2em
\Phi^i_{\!j}(-\varepsilon)=\delta^i_j-\frac{\partial X^i(x^0,x^1,x^2,
x^3)}{\partial x^j}\ \varepsilon+\,\ldots.
\mytag{2.13}
\endalign
$$
Here $X^i(x^0,x^1,x^2,x^3)$ are the components of the vector field 
$\bold X$ in the local coordinates $x^0,\,x^1,\,x^2,\,x^3$. Now if we
denote by $L_{\bold X}(\bold Y)^{i_1\ldots\,i_r}_{j_1\ldots\,j_s}$ the
components of $L_{\bold X}(\bold Y)$, then, applying \mythetag{2.12} and 
\mythetag{2.12} to \mythetag{2.8} and taking into account \mythetag{2.9}, 
we obtain
$$
\hskip -2em
\gathered
L_{\bold X}(\bold Y)^{i_1\ldots\,i_r}_{j_1\ldots\,j_s}
=\sum^s_{m=1}
\sum^3_{k_m=0}\frac{\partial X^{k_m}}{\partial x^{j_m}}
\ Y^{\,i_1\ldots\,\ldots\,\ldots\,
i_r}_{j_1\ldots\,k_m\ldots\,j_s}-\\
-\sum^r_{m=1}\sum^3_{k_m=0}\frac{\partial X^{i_m}}{\partial x^{k_m}}
\ Y^{\,i_1\ldots\,k_m\ldots\,
i_r}_{j_1\ldots\,\ldots\,\ldots\,j_s}
+\sum^3_{k=0}X^k
\,\frac{\partial Y^{i_1\ldots\,i_r}_{j_1\ldots\,j_s}}{\partial x^k}.
\endgathered
\mytag{2.14}
$$
The Lie derivative $L_{\bold X}$ given by the formula \mythetag{2.14}
possesses the following properties:
\roster
\rosteritemwd=5pt
\item the Lie derivative $L_{\bold X}$ preserves the type of a 
      tensor field, i\.\,e\. $\bold Y$ and $L_{\bold X}(\bold Y)$ 
      are tensor fields of the same type;
\item $L_{\bold X}(\bold Y_1\otimes\bold Y_2)=L_{\bold X}(\bold Y_1)
      \otimes\bold Y_2+\bold Y_1\otimes L_{\bold X}(\bold Y_2)$ 
      for arbitrary two tensorial fields $\bold Y_1$ and $\bold Y_2$;
\item $L_{\bold X}(C(\bold Y))=C(L_{\bold X}(\bold Y))$, i\.\,e\.
      $L_{\bold X}$ commute with contractions.
\endroster
The properties \therosteritem{1}--\therosteritem{3} are easily derived 
with the use of the formula \mythetag{2.14} itself.
\head
3. Lie derivatives in frames formalism.
\endhead
     Let $x^0,\,x^1,\,x^2,\,x^3$ be the local coordinates of some local 
chart of the space-time manifold $M$. The coordinates $x^0,\,x^1,\,x^2,
\,x^3$ induce the frame $\bold X_0,\,\bold X_1,\,\bold X_2,\,\bold X_3$ 
of the coordinate vector fields in the domain of these coordinates:
$$
\xalignat 4
&\hskip -2em
\bold X_0=\frac{\partial}{\partial x^0},
&&\bold X_1=\frac{\partial}{\partial x^1},
&&\bold X_2=\frac{\partial}{\partial x^2},
&&\bold X_3=\frac{\partial}{\partial x^3}.
\qquad
\mytag{3.1}
\endxalignat
$$
The frame composed by the vector fields \mythetag{3.1} is a holonomic
frame since these vector fields commute with each other:
$$
[\bold X_i,\,\bold X_j]=0.
$$
However, one can consider some non-holonomic frame 
$\boldsymbol\Upsilon_0,\,\boldsymbol\Upsilon_1,\,\boldsymbol
\Upsilon_2,\,\boldsymbol\Upsilon_3$, i\.\,e\. a frame with
non-commuting vector fields:
$$
\hskip -2em
[\boldsymbol\Upsilon_{\!i},\boldsymbol\Upsilon_{\!j}]=\sum^3_{k=0}
c^{\,k}_{ij}\,\boldsymbol\Upsilon_{\!k}.
\mytag{3.2}
$$
The commutation coefficients $c^{\,k}_{ij}$ in \mythetag{3.2} are uniquely 
determined by the frame vector fields $\boldsymbol\Upsilon_0,\,\boldsymbol
\Upsilon_1,\,\boldsymbol\Upsilon_2,\,\boldsymbol\Upsilon_3$ since they are
linearly independent at each point of their domain. Our nearest goal is to
derive the formula analogous to \mythetag{2.14} for the case where all
tensorial fields are represented by their components is some non-holonomic
frame $\boldsymbol\Upsilon_0,\,\boldsymbol\Upsilon_1,\,\boldsymbol\Upsilon_2,
\,\boldsymbol\Upsilon_3$.\par
     Let $\varphi$ be a scalar field. Then the Lie derivative $L_{\bold X}$
of $\varphi$ is reduced to the differentiation of the function $\varphi$
along the vector $\bold X$:
$$
\hskip -2em
L_{\bold X}(\varphi)=\sum^3_{k=0}X^k\,\frac{\partial\varphi}{\partial x^k}.
\mytag{3.3}
$$
The formula \mythetag{3.3} is easily derived by substituting $\bold Y=\varphi$
with $r=0$ and $s=0$ into \mythetag{2.14}. Similarly, if $\bold Y$ is a vector
field, from \mythetag{2.14} we derive
$$
L_{\bold Y}(\bold Y)=[\bold X,\bold Y].
\mytag{3.4}
$$
Now assume that both of the vector fields $\bold X$ and $\bold Y$ are 
represented by their expansions in a non-holonomic frame $\boldsymbol
\Upsilon_0,\,\boldsymbol\Upsilon_1,\,\boldsymbol\Upsilon_2,\,\boldsymbol
\Upsilon_3$:
$$
\xalignat 2
&\hskip -2em
\bold X=\sum^3_{k=0}X^k\,\,\boldsymbol\Upsilon_{\!k},
&&\bold Y=\sum^3_{i=0}Y^i\,\,\boldsymbol\Upsilon_{\!i}.
\mytag{3.5}
\endxalignat
$$
Then, substituting \mythetag{3.5} into \mythetag{3.4}, we derive
$$
L_{\bold X}(\bold Y)^k=\sum^3_{i=0}X^i\,L_{\boldsymbol\Upsilon_i}(Y^k)
-\sum^3_{i=0}Y^i\,L_{\boldsymbol\Upsilon_i}(X^k)+\sum^3_{i=0}\sum^3_{j=0}
c^{\,k}_{ij}\,X^i\,Y^j.
\quad
\mytag{3.6}
$$
The Lie derivatives $L_{\boldsymbol\Upsilon_i}(Y^k)$ and 
$L_{\boldsymbol\Upsilon_i}(X^k)$ in \mythetag{3.6} are calculated
according to the formula \mythetag{3.3}, i\.\,e\. we substitute
$\varphi=Y^k$ and $\varphi=X^k$ into \mythetag{3.3}. Note that $X^k$ 
in \mythetag{3.3} differ from that of \mythetag{3.5}. The components 
of $\bold X$ in the formula \mythetag{3.3} are taken from the expansion 
of the vector field $\bold X$ in the holonomic frame \mythetag{3.1}.
The non-holonomic version of the formula  \mythetag{3.3} looks like
$$
L_{\bold X}(\varphi)=\sum^3_{k=0}X^k\,L_{\boldsymbol\Upsilon_i}
(\varphi).
$$\par
    Let's denote by $\boldsymbol\eta^0,\,\boldsymbol\eta^1,\,\boldsymbol
\eta^2,\,\boldsymbol\eta^3$ the dual frame for $\boldsymbol\Upsilon_0,\,
\boldsymbol\Upsilon_1,\,\boldsymbol\Upsilon_2,\,\boldsymbol\Upsilon_3$.
This means that $\boldsymbol\eta^0,\,\boldsymbol\eta^1,\,\boldsymbol\eta^2,
\,\boldsymbol\eta^3$ are four covectorial fields such that 
$$
\hskip -2em
\eta^i(\boldsymbol\Upsilon_{\!j})=\bigl(\boldsymbol\eta^i,\,\boldsymbol
\Upsilon_{\!j}\bigr)=C(\boldsymbol\eta^i\otimes\boldsymbol\Upsilon_{\!j})
=\delta^i_j.
\mytag{3.7}
$$
Using the properties \therosteritem{1}--\therosteritem{3} from Section~2
and using \mythetag{3.4}, from \mythetag{3.7} we derive 
$$
\hskip -2em
L_{\boldsymbol\Upsilon_i}(\boldsymbol\eta^k)=-\sum^3_{j=0}c^{\,k}_{ij}
\,\boldsymbol\eta^j.
\mytag{3.8}
$$
The commutation coefficients $c^{\,k}_{ij}$ in \mythetag{3.8} are the 
same as in \mythetag{3.2}. Assume that $\bold Y$ is a covectorial field
expanded in the frame $\boldsymbol\eta^0,\,\boldsymbol\eta^1,\,\boldsymbol
\eta^2,\,\boldsymbol\eta^3$:
$$
\hskip -2em
\bold Y=\sum^3_{j=1}Y_i\,\boldsymbol\eta^i.
\mytag{3.9}
$$
From \mythetag{3.7}, \mythetag{3.8}, and \mythetag{3.9} we derive
$$
L_{\bold X}(\bold Y)_k=\sum^3_{i=0}X^i\,L_{\boldsymbol\Upsilon_i}(Y_k)
+\sum^3_{i=0}Y_i\,L_{\boldsymbol\Upsilon_k}(X^i)-\sum^3_{i=0}\sum^3_{j=0}
X^i\,Y_j\,c^{\,j}_{ik}.
\quad
\mytag{3.10}
$$ 
Now, using the properties \therosteritem{1}--\therosteritem{3} again
and combining \mythetag{3.6} with \mythetag{3.10}, we can extend the formula 
\mythetag{3.6} to the case of an arbitrary tensor field $\bold Y$ of the type 
$(r,s)$:
$$
\hskip -2em
\gathered
L_{\bold X}(\bold Y)^{i_1\ldots\,i_r}_{j_1\ldots\,j_s}
=\sum^3_{i=0}X^i\,L_{\boldsymbol\Upsilon_i}(Y^{i_1\ldots
\,i_r}_{j_1\ldots\,j_s})\,+\\
+\sum^s_{m=1}\sum^3_{k_m=0}
\!\left(\!L_{\boldsymbol\Upsilon_{j_m}}(X^{k_m})
-\sum^3_{i=0}c^{\,k_m}_{ij_m}\,X^i\!\right)Y^{\,i_1\ldots\,\ldots\,\ldots\,
i_r}_{j_1\ldots\,k_m\ldots\,j_s}\,-\\
-\sum^r_{m=1}\sum^3_{k_m=0}
\!\left(\!L_{\boldsymbol\Upsilon_{k_m}}(X^{i_m})
-\sum^3_{i=0}c^{\,i_m}_{ik_m}\,X^i
\!\right)Y^{\,i_1\ldots\,k_m\ldots\,
i_r}_{j_1\ldots\,\ldots\,\ldots\,j_s}.
\endgathered
\mytag{3.11}
$$
The formulas \mythetag{3.6} and \mythetag{3.10} are special cases of the
formula \mythetag{3.11}. If the frame $\boldsymbol\Upsilon_0,\,\boldsymbol
\Upsilon_1,\,\boldsymbol\Upsilon_2,\,\boldsymbol\Upsilon_3$ coincides with
the holonomic frame \mythetag{3.1}, then $c^{\,k}_{ij}=0$ and the formula 
\mythetag{3.11} reduces to \mythetag{2.14}.\par
      Now let's return back to the formula \mythetag{2.8}. This formula is
valid in frame presentation of tensor fields too. However, the matrices 
$\Phi^i_{\!j}(\varepsilon)$ and $\Phi^i_{\!j}(-\varepsilon)$ in this
case are not given by the formulas \mythetag{2.7}. Here we use the
formulas 
$$
\align
&\hskip -2em
\Phi^i_{\!j}(\varepsilon)=\delta^i_j
+\!\left(\!L_{\boldsymbol\Upsilon_j}(X^i)-\sum^3_{m=0}X^m\,c^{\,i}_{mj}
\!\right)\varepsilon+\,\ldots,
\mytag{3.12}\\
&\hskip -2em
\Phi^i_{\!j}(-\varepsilon)=\delta^i_j
-\!\left(\!L_{\boldsymbol\Upsilon_j}(X^i)-\sum^3_{m=0}X^m\,c^{\,i}_{mj}
\!\right)\varepsilon+\,\ldots.
\mytag{3.13}
\endalign
$$
The formulas \mythetag{3.12} and \mythetag{3.13} are analogous to
\mythetag{2.12} and \mythetag{2.13}.
\head
4. Tangent vector fields on vector bundles.
\endhead
      Let $VM$ be an $n$-dimensional vector bundle over the space-time 
manifold $M$. Any trivialization of $VM$ is given by $n$ sections 
$\boldsymbol\Upsilon_{\!1},\,\ldots,\,\boldsymbol\Upsilon_{\!n}$ linearly
independent at each point of their domain $U$. Let $\bold v$ be a vector
of the fiber $V_q(M)$:
$$
\hskip -2em
\bold v=v^1\,\boldsymbol\Upsilon_{\!1}+\ldots
+v^n\,\boldsymbol\Upsilon_{\!n}.
\mytag{4.1}
$$
Then $\tilde q=(q,\bold v)$ is a point of $VM$. If $x^0,\,x^1,\,x^2,\,x^3$ 
are some local coordinates within the domain $U\subset M$ and $v^0,\,v^1,
\,v^2,\,v^3$ are taken from \mythetag{4.1}, then 
$$
\hskip -2em
x^0,\,\ldots,\,x^3,\,v^1,\,\ldots,\,v^n
\mytag{4.2}
$$
are the coordinates of the point $\tilde q=(q,\bold v)$. The coordinates 
\mythetag{4.2} are naturally subdivided into two groups --- the base
coordinates $x^0,\,\ldots,\,x^3$ and the fiber coordinates $v^1,\,\ldots,
\,v^n$. Let $\bundleboldX$ be a tangent vector field on $VM$. In the local 
coordinates \mythetag{4.2} it is represented as the following differential 
operator:
$$
\hskip -2em
\bundleboldX=\sum^3_{i=0}X^i\,\frac{\partial}{\partial x^i}
+\sum^n_{i=1}V^i\,\frac{\partial}{\partial v^i}.
\mytag{4.3}
$$
Under the canonical projection $\pi\!:\,VM\to M$ the vector field 
\mythetag{4.3} is mapped to 
$$
\hskip -2em
\bold X=\sum^3_{i=0}X^i\,\frac{\partial}{\partial x^i}.
\mytag{4.4}
$$
The vector field \mythetag{4.3} produces the one-parametric 
local group of diffeomorphisms
$$
\hskip -2em
\varphi_{\varepsilon}\!:\,VM\to VM
\mytag{4.5}
$$
that extends the local group \mythetag{2.1} produced by the vector field
\mythetag{4.4}. Due to \mythetag{4.5} the functions \mythetag{2.6} are 
complemented with the functions
$$
\align
&\hskip -2em
\left\{\aligned
&v^1=U^1(\varepsilon,y^0,\,\ldots,y^3,w^1,\,\ldots,w^n),\\
\vspace{-1ex}
&.\ .\ .\ .\ .\ .\ .\ .\ .\ .\ .\ .\ .\ .\ .\ .\ .\ .\ 
.\ .\ .\ .\ .\ .\ .\ \\
\vspace{-0.5ex}
&v^n=U^n(\varepsilon,y^0,\,\ldots,y^3,w^1,\,\ldots,w^n),
\endaligned\right.
\mytag{4.6}\\
\vspace{2ex}
&\hskip -2em
\left\{\aligned
&w^1=U^1(-\varepsilon,x^0,\,\ldots,x^3,v^1,\,\ldots,v^n),\\
&.\ .\ .\ .\ .\ .\ .\ .\ .\ .\ .\ .\ .\ .\ .\ .\ .\ .\ 
.\ .\ .\ .\ .\ .\ .\ \\
&w^n=U^n(-\varepsilon,x^0,\,\ldots,x^3,v^1,\,\ldots,v^n).
\endaligned\right.
\mytag{4.7}
\endalign
$$
\mydefinition{4.1} The tangent vector field \mythetag{4.3} on 
$VM$ is called {\it concordant with the bundle structure\/} if 
the functions \mythetag{4.6} and \mythetag{4.7} are linear with 
respect to their arguments $v^1,\,\ldots,\,v^n$ and $w^1,\,\ldots,\,w^n$.
\enddefinition
     In the case of a concordant vector field $\bundleboldX$ the local
diffeomorphisms \mythetag{4.5} break into the series of linear mappings
$$
\hskip -2em
\varphi_{\varepsilon}\!:\,V_p(M)\to V_q(M),
\mytag{4.8}
$$
where $q=\varphi_{\varepsilon}(p)$. The functions \mythetag{4.6}
present the diffeomorphisms \mythetag{4.8} in local coordinates. 
In the case of a concordant vector field they are given by the 
formulas
$$
\hskip -2em
U^i(\varepsilon,y^0,\,\ldots,y^3,w^1,\,\ldots,w^n)
=\sum^n_{j=1}U^i_j(\varepsilon,y^0,\,\ldots,y^3)\,w^j.
\mytag{4.9}
$$
The vertical components of the tangent vector field \mythetag{4.3}
are given by the derivatives
$$
\hskip -2em
V^i=\frac{dU^i(\varepsilon,x^0,\,\ldots,x^3,v^1,\,\ldots,v^n)}
{d\varepsilon}\,\hbox{\vrule height 12pt depth 8pt width 0.5pt}_{\ 
\varepsilon=0}.
\mytag{4.10}
$$
Substituting \mythetag{4.9} into \mythetag{4.10}, we derive 
$$
\hskip -2em
V^i=\sum^n_{j=1}V^i_j(x^0,\,\ldots,x^3)\,v^j\text{, \ where \ }
V^i_j=U^i_{\!j}\,\hbox{\vrule height 8pt depth 8pt width 0.5pt}_{\ 
\varepsilon=0}.
\mytag{4.11}
$$
\mytheorem{4.1} The tangent vector field \mythetag{4.3} on a vector
bundle $VM$ is concordant with the bundle structure if and only if
its vertical components are linear functions with respect to
$v^1,\,\ldots,\,v^n$ given by the formula \mythetag{4.11}
\endproclaim
     For the matrices $U^i_{\!j}$ in \mythetag{4.9}, which represent
the linear mappings \mythetag{4.8} in the frame $\boldsymbol
\Upsilon_{\!1},\,\ldots,\,\boldsymbol\Upsilon_{\!n}$, the formulas
\mythetag{4.10} and \mythetag{4.11} yield
$$
\hskip -2em
\aligned
&U^i_{\!j}(\varepsilon)=\delta^i_j+V^i_j(x^0,\,\ldots,x^3)\,\varepsilon
+\,\ldots,\\
\vspace{2ex}
&U^i_{\!j}(-\varepsilon)=\delta^i_j-V^i_j(x^0,\,\ldots,x^3)\,\varepsilon
+\,\ldots.
\endaligned
\mytag{4.12}
$$
The expansions \mythetag{4.12} are similar to \mythetag{2.12}, 
\mythetag{2.13}, \mythetag{3.12}, and \mythetag{3.13}.\par
     Let $\bold Y$ be a tensor field of the type $(r,s)$ associated 
with the vector bundle $VM$, i\.\,e\. let $\bold Y$ be a section of 
$V^r_sM$, where $V^r_sM$ is the following tensor bundle:
$$
\hskip -2em
V^r_sM=\overbrace{VM\otimes\ldots\otimes VM}^{\text{$r$ times}}
\otimes\underbrace{V^*\!M\otimes\ldots\otimes V^*\!M}_{\text{$s$ times}}.
\mytag{4.13}
$$
Assume that $Y^{i_1\ldots\,i_r}_{j_1\ldots\,j_s}(x^0,x^1,x^2,x^3)$ are
the components of the tensor field $\bold Y$ in the frame $\boldsymbol
\Upsilon_{\!1},\,\ldots,\,\boldsymbol\Upsilon_{\!n}$. Then, using the
quantities $V^i_j$ from the expansions \mythetag{4.12}, we define the 
Lie derivative $\Cal L_{\bold X}(\bold Y)$ of the field $\bold Y$:
$$
\hskip -2em
\gathered
\Cal L_{\bold X}(\bold Y)^{i_1\ldots\,i_r}_{j_1\ldots\,j_s}
=\sum^3_{i=0}X^i\,L_{\boldsymbol\Upsilon_i}(Y^{i_1\ldots
\,i_r}_{j_1\ldots\,j_s})\,+\\
+\sum^s_{m=1}\sum^n_{k_m=1}V^{k_m}_{j_m}
\,Y^{\,i_1\ldots\,\ldots\,\ldots\,
i_r}_{j_1\ldots\,k_m\ldots\,j_s}
-\sum^r_{m=1}\sum^n_{k_m=1}V^{i_m}_{k_m}
\,Y^{\,i_1\ldots\,k_m\ldots\,
i_r}_{j_1\ldots\,\ldots\,\ldots\,j_s}.
\endgathered
\mytag{4.14}
$$
The Lie derivative \mythetag{4.14} is called natural if the quantities
$V^i_j$ are expressed in some natural way through the components of the
vector field $\bold X$ in \mythetag{4.4}.
\head
5. Natural liftings and Kosmann liftings.
\endhead
     In this section we apply the results of the previous section~4
to the tangent bundle $TM$, i\.\,e\. we set $VM=TM$. Comparing
the formula \mythetag{4.13} with \mythetag{2.4} and the formula 
\mythetag{4.14} with \mythetag{3.11}, we find that
$$
\hskip -2em
V^i_j=L_{\boldsymbol\Upsilon_j}(X^i)-\sum^3_{m=0}X^m\,c^{\,i}_{mj}.
\mytag{5.1}
$$
Now we substitute \mythetag{5.1} into \mythetag{4.11} and then we
substitute \mythetag{4.11} into \mythetag{4.3}: 
$$
\bundleboldX_N=\sum^3_{i=0}\sum^3_{m=0}X^m\,\varUpsilon^i_m
\,\frac{\partial}{\partial x^i}
+\sum^3_{i=0}\sum^3_{j=0}\left(\!L_{\boldsymbol\Upsilon_j}(X^i)
-\sum^3_{m=0}X^m\,c^{\,i}_{mj}\right)v^j\,
\frac{\partial}{\partial v^i}.\quad
\mytag{5.2}
$$
The tangent vector field \mythetag{5.2} on $TM$ is a natural lifting
of the vector field
$$
\hskip -2em
\bold X=\sum^3_{i=0}\sum^3_{m=0}X^m\,\varUpsilon^i_m
\,\frac{\partial}{\partial x^i}=\sum^3_{m=0}X^m\,
\boldsymbol\Upsilon_{\!m}
\mytag{5.3}
$$
from $M$ to $TM$. The Lie derivative \mythetag{4.14} determined by
the vector field \mythetag{5.3} and by its lifting \mythetag{5.2} 
is a natural Lie derivative coinciding with \mythetag{3.11} and 
\mythetag{2.14}.\par
      Now let's recall that the the space-time manifold $M$ is equipped
with the metric $\bold g$. Its signature is $(+,-,-,-)$. For this reason
each fiber $T_p(M)$ of the tangent bundle $TM$ is a pseudo-Euclidean 
linear vector space.
\mydefinition{5.1} A lifting $\bundleboldX$ of a vector field $\bold X$
from $M$ to $TM$ is called a {\it Kosmann lifting\/} if the linear 
mappings \mythetag{4.8} associated with this lifting are isometries.
\enddefinition
     Kossman liftings were first introduced by Yvette Kosmann in
\myciterange{4}{4}{--}{7}.
\mytheorem{5.1} The natural lifting \mythetag{5.2} of a vector field\/
$\bold X$ is a Kosmann lifting if and only if\/ $\bold X$ is a Killing 
vector field.
\endproclaim
     The proof is trivial. By definition, Killing vector fields are those 
whose local diffeomorphisms preserve the metric tensor $\bold g$. Hence,
the mapping $\varphi_{\varepsilon*}$ from \mythetag{2.2} restricted
to any fiber $T_p(M)$ is a linear isometry.\par
     Let's study the isometry condition from the 
definition~\mythedefinition{5.1} in more details. Applying the linear 
mappings \mythetag{4.8} to the metric tensor, we get the equality
$$
\hskip -2em
g_{ij}(x^0,\,\ldots,x^3)=\sum^3_{r=0}\sum^3_{s=0}
U^r_i(-\varepsilon)\ U^s_{\!j}(-\varepsilon)\ g_{rs}(y^0,\,\ldots,y^3).
\mytag{5.4}
$$
Differentiating the formula \mythetag{5.4} with respect to $\varepsilon$, 
we take into account the formulas \mythetag{4.12}, \mythetag{2.6}, 
\mythetag{2.10}, \mythetag{2.11}, \mythetag{4.10}, and \mythetag{4.11}.
As a result we get
$$
\hskip -2em
0=\sum^3_{r=0}V^r_i\,g_{rj}+\sum^3_{r=0}V^r_j\,g_{ir}
+\sum^3_{m=0}X^m\,L_{\boldsymbol\Upsilon_m}(g_{ij}).
\mytag{5.5}
$$
The equality \mythetag{5.5} can be simplified to
$$
\hskip -2em
\Cal L_{\bold X}(\bold g)=0,
\mytag{5.6}
$$
where the Lie derivative $\Cal L_{\bold X}$ is calculated according to 
the formula \mythetag{4.14}. We shall treat the equality \mythetag{5.6}
neither as a condition for $\bold X$ nor as a condition for $\bold g$,
but as a condition for $V^i_j$. For this purpose we denote
$$
\hskip -2em
V_{ij}=\sum^3_{r=0}V^r_i\,g_{rj}.
\mytag{5.7}
$$
Then the equality \mythetag{5.5}, which is  equivalent to \mythetag{5.6}
is written as follows:
$$
\hskip -2em
V_{ij}+V_{j\kern 0.5pt i}=-\sum^3_{m=0}X^m\,L_{\boldsymbol
\Upsilon_m}(g_{ij}).
\mytag{5.8}
$$
The equality  \mythetag{5.8} fixes the symmetric part of $V_{ij}$ for 
a Kosmann lifting of a vector field. The skew-symmetric part of $V_{ij}$
can be obtained by alternating \mythetag{5.1}. As a result we get the
following two formulas:
$$
\align
&\hskip -2em
V^{\sssize\text{sym}}_{ij}=-\frac{1}{2}\sum^3_{m=0}X^m\,L_{\boldsymbol
\Upsilon_m}(g_{ij}),
\mytag{5.9}\\
&\hskip -2em
\aligned
V^{\sssize\text{skew}}_{ij}=&\frac{1}{2}\sum^3_{r=0}L_{\boldsymbol
\Upsilon_i}(X^r)\,g_{rj}-\frac{1}{2}\sum^3_{r=0}L_{\boldsymbol\Upsilon_j}
(X^r)\,g_{ri}\,-\\
&-\frac{1}{2}\sum^3_{r=0}\sum^3_{m=0}X^m\,c^{\,r}_{mi}\,g_{rj}
+\frac{1}{2}\sum^3_{r=0}\sum^3_{m=0}X^m\,c^{\,r}_{mj}\,g_{ri}.
\endaligned
\mytag{5.10}
\endalign
$$
Adding the formulas \mythetag{5.9} and \mythetag{5.10}, we derive
$$
\gathered
V_{ij}=-\frac{1}{2}\sum^3_{m=0}X^m\,L_{\boldsymbol
\Upsilon_m}(g_{ij})+\frac{1}{2}\sum^3_{r=0}L_{\boldsymbol
\Upsilon_i}(X^r)\,g_{rj}\,-\\
-\,\frac{1}{2}\sum^3_{r=0}L_{\boldsymbol\Upsilon_j}(X^r)\,g_{ri}
-\frac{1}{2}\sum^3_{r=0}\sum^3_{m=0}X^m\,c^{\,r}_{mi}\ g_{rj}
+\frac{1}{2}\sum^3_{r=0}\sum^3_{m=0}X^m\,c^{\,r}_{mj}\ g_{ri}.
\endgathered
\qquad
\mytag{5.11}
$$
In order to get back to $V^j_{\,i}$ we need to raise the index $j$ in 
\mythetag{5.11}:
$$
\gathered
V^j_{\,i}=-\frac{1}{2}\sum^3_{m=0}\sum^3_{r=0}X^m
\,L_{\boldsymbol\Upsilon_m}(g_{ir})\,g^{rj}
-\frac{1}{2}\sum^3_{r=0}\sum^3_{s=0}L_{\boldsymbol\Upsilon_s}(X^r)
\,g_{ri}\,g^{sj}\,+\\
+\,\frac{1}{2}\,L_{\boldsymbol\Upsilon_i}(X^j)-\frac{1}{2}\sum^3_{m=0}
X^m\,c^{\,j}_{mi}+\frac{1}{2}\sum^3_{r=0}\sum^3_{s=0}\sum^3_{m=0}
X^m\,c^{\,r}_{ms}\,g_{ri}\,g^{sj}.
\endgathered
\qquad
\mytag{5.12}
$$
Note that in order to fit \mythetag{4.11} we should exchange the
indices $i$ and $j$ in \mythetag{5.12}:
$$
\gathered
V^i_j=-\frac{1}{2}\sum^3_{m=0}\sum^3_{r=0}g^{ir}\,X^m
\,L_{\boldsymbol\Upsilon_m}(g_{rj})-\frac{1}{2}\sum^3_{r=0}
\sum^3_{s=0}g^{is}\,L_{\boldsymbol\Upsilon_s}(X^r)\,g_{rj}\,+\\
+\,\frac{1}{2}\,L_{\boldsymbol\Upsilon_j}(X^i)-\frac{1}{2}\sum^3_{m=0}
X^m\,c^{\,i}_{mj}+\frac{1}{2}\sum^3_{r=0}\sum^3_{s=0}\sum^3_{m=0}
g^{is}\,X^m\,c^{\,r}_{ms}\,g_{rj}.
\endgathered
\qquad
\mytag{5.13}
$$
As a result, substituting \mythetag{5.13} into \mythetag{4.11}, we derive 
the formula
$$
\gathered
V^i=-\frac{1}{2}\sum^3_{m=0}\sum^3_{r=0}\sum^3_{j=0}g^{ir}\,X^m
\,L_{\boldsymbol\Upsilon_m}(g_{rj})\,v^j
+\frac{1}{2}\sum^3_{j=0}L_{\boldsymbol\Upsilon_j}(X^i)\,v^j\,-\\
-\,\frac{1}{2}\sum^3_{r=0}\sum^3_{s=0}\sum^3_{j=0}g^{is}
\,L_{\boldsymbol\Upsilon_s}(X^r)\,g_{rj}\,v^j
-\frac{1}{2}\sum^3_{m=0}\sum^3_{j=0}X^m\,c^{\,i}_{mj}\,v^j\,+\\
+\,\frac{1}{2}\sum^3_{r=0}\sum^3_{s=0}\sum^3_{m=0}\sum^3_{j=0}
g^{is}\,X^m\,c^{\,r}_{ms}\,g_{rj}\,v^j.
\endgathered
\qquad
\mytag{5.14}
$$
Applying \mythetag{5.14} to \mythetag{4.3} and taking into account that
we choose $VM=TM$, we get the following tangent vector field on the
tangent bundle:
$$
\hskip -0.1em
\gathered
\bundleboldX_K=\sum^3_{i=0}\sum^3_{m=0}X^m\,\varUpsilon^i_m
\,\frac{\partial}{\partial x^i}
+\frac{1}{2}\sum^3_{i=0}\left(\,\sum^3_{j=0}L_{\boldsymbol\Upsilon_j}
(X^i)\,v^j\,-\right.\\
\left.-\sum^3_{m=0}\sum^3_{j=0}X^m\,c^{\,i}_{mj}\,v^j
+\sum^3_{r=0}\sum^3_{s=0}\sum^3_{m=0}\sum^3_{j=0}
g^{is}\,X^m\,c^{\,r}_{ms}\,g_{rj}\,v^j\,-\right.\\
\left.-\sum^3_{r=0}\sum^3_{s=0}\sum^3_{j=0}g^{is}
\,L_{\boldsymbol\Upsilon_s}(X^r)\,g_{rj}\,v^j
-\sum^3_{m=0}\sum^3_{r=0}\sum^3_{j=0}g^{ir}\,X^m
\,L_{\boldsymbol\Upsilon_m}(g_{rj})\,v^j\!\right)
\frac{\partial}{\partial v^i}.
\endgathered
\hskip -1em
\mytag{5.15}
$$
\mydefinition{5.2} The tangent vector field \mythetag{5.15} on the
tangent bundle $TM$ is called the {\it standard Kosmann lifting\/}
of the vector field \mythetag{4.4} from $M$ to $TM$.
\enddefinition
     We can calculate the standard Kosmann lifting \mythetag{5.15}
using the holonomic frame \mythetag{3.1} instead of the non-holonomic
frame $\boldsymbol\Upsilon_0,\,\boldsymbol\Upsilon_1,\,\boldsymbol
\Upsilon_2,\,\boldsymbol\Upsilon_3$. Then \mythetag{5.15} reduces to
$$
\hskip -2em
\gathered
\bundleboldX_K=\sum^3_{i=0}X^i\,\frac{\partial}{\partial x^i}
+\frac{1}{2}\sum^3_{i=0}\!\left(-\sum^3_{r=0}\sum^3_{s=0}
\sum^3_{j=0}g^{is}
\,\frac{\partial X^r}{\partial x^s}\,g_{rj}\ v^j\,+\right.\\
\left.+\sum^3_{j=0}\frac{\partial X^i}{\partial x^j}\,v^j
-\sum^3_{m=0}\sum^3_{r=0}\sum^3_{j=0}g^{ir}\,X^m
\,\frac{\partial g_{rj}}{\partial x^m}\,v^j\!\right)
\!\frac{\partial}{\partial v^i}.
\endgathered
\mytag{5.16}
$$\par
     Now let's remember that the metric $\bold g$ is associated with
the metric connection $\Gamma$. It is known as the Levi-Civita
connections. The components of the Levi-Civita connection in a holonomic 
frame \mythetag{3.1} are given by the formula
$$
\hskip -2em
\Gamma^k_{ij}=\sum^3_{r=0}\frac{g^{kr}}{2}\left(
\frac{\partial g_{rj}}{\partial x^i}+\frac{\partial g_{ir}}{\partial x^j}
-\frac{\partial g_{ij}}{\partial x^r}\right).
\mytag{5.17}
$$
The formula \mythetag{5.17} is easily derived from the following two 
conditions (see \mycite{8}):
$$
\xalignat 2
&\hskip -2em
\Gamma^k_{ij}=\Gamma^k_{ji},
&&\nabla_{\!k}g_{ij}=0.
\mytag{5.18}
\endxalignat
$$
Using the connection components $\Gamma^k_{ij}$, now we express the 
partial derivatives in the formula \mythetag{5.16} through the 
corresponding covariant derivatives:
$$
\hskip -2em
\aligned
&\frac{\partial X^r}{\partial x^s}=\nabla_{\!s}X^r-
\sum^3_{m=0}\Gamma^r_{sm}\,X^m,\\
&\frac{\partial X^i}{\partial x^j}=\nabla_{\!j}X^i-
\sum^3_{m=0}\Gamma^i_{jm}\,X^m.
\endaligned
\mytag{5.19}
$$
Moreover, from the second equality \mythetag{5.18} we derive
$$
\hskip -2em
\frac{\partial g_{rj}}{\partial x^m}=
\sum^3_{s=0}\Gamma^s_{mj}\,g_{rs}
+\sum^3_{s=0}\Gamma^s_{mr}\,g_{sj}.
\mytag{5.20}
$$
Substituting \mythetag{5.19} and \mythetag{5.20} back into 
\mythetag{5.16} and taking into account the first equality
in \mythetag{5.17}, we get the following formula:
$$
\hskip -2em
\gathered
\bundleboldX_K=\sum^3_{i=0}X^i\,\frac{\partial}{\partial x^i}
+\sum^3_{i=0}\sum^3_{j=0}\!\left(-\frac{1}{2}\sum^3_{r=0}
\sum^3_{s=0}g^{is}\,\nabla_{\!s}X^r\,g_{rj}\,+\right.\\
\left.+\,\frac{1}{2}\nabla_{\!j}X^i
-\sum^3_{m=0}X^m\,\Gamma^i_{mj}\!\right)\!v^j
\,\frac{\partial}{\partial v^i}.
\endgathered
\mytag{5.21}
$$\par
     The components of the natural lifting \mythetag{5.2} also can
be expressed through the connection components $\Gamma^k_{ij}$ and
covariant derivatives in the holonomic frame \mythetag{3.1}:
$$
\hskip -2em
\gathered
\bundleboldX_N=\sum^3_{i=0}X^i\,\frac{\partial}{\partial x^i}
+\sum^3_{i=0}\sum^3_{j=0}\!\left(\!\nabla_{\!j}X^i
-\sum^3_{m=0}X^m\,\Gamma^i_{mj}\!\right)\!v^j\,\frac{\partial}
{\partial v^i}.
\endgathered
\mytag{5.22}
$$
Two different liftings \mythetag{5.21} and \mythetag{5.22} are 
associated with two different Lie derivatives $\Cal L_{\bold X}$
and $L_{\bold X}$ respectively. Here $L_{\bold X}$ is the regular
Lie derivative, while $\Cal L_{\bold X}$ is called the {\it standard
Kosmann-Lie derivative}. Both of them are differentiations of the
algebra of tensor fields. Comparing \mythetag{5.21} and \mythetag{5.22},
we get the relationship
$$
\hskip -2em
\Cal L_{\bold X}=L_{\bold X}+S_{\bold X}.
\mytag{5.23}
$$
Here $S_{\bold X}$ is a degenerate differentiation in the sense of 
the proposition~3.3 in Chapter~\uppercase\expandafter{\romannumeral 1}
of \mycite{3}. From \mythetag{5.21} and \mythetag{5.22} we derive
that the degenerate differentiation $S_{\bold X}$ in \mythetag{5.23} 
is given by the tensor field $\bold S_{\bold X}$ with the following 
components:
$$
\hskip -2em
S^{\kern 0.2pt i}_j({\bold X})=\frac{\nabla_{\!j}X^i+\nabla^iX_j}{2}
=\frac{1}{2}\,\nabla_{\!j}X^i+\frac{1}{2}\sum^3_{r=0}
\sum^3_{s=0}g^{is}\,\nabla_{\!s}X^r\,g_{rj}
\mytag{5.24}
$$
The tensor field $\bold S_{\bold X}$ with the components
\mythetag{5.24} is equal to zero if and only if\/ $\bold X$ is 
a Killing vector field. This fact is easily derived from the
theorem~\mythetheorem{5.1} or from the equality \mythetag{5.6},
which is fulfilled identically by definition.\par
     Another special case for the formula \mythetag{5.15} is the
case of a non-holonomic, but orthonormal frame $\boldsymbol
\Upsilon_0,\,\boldsymbol\Upsilon_1,\,\boldsymbol\Upsilon_2,
\,\boldsymbol\Upsilon_3$. In such a frame the components of
the metric tensor are constants. They are given by the 
Minkowski matrix:
$$
\hskip -2em
g_{ij}=g^{ij}=\Vmatrix 1 & 0 & 0 & 0\\0 & -1 & 0 & 0\\
0 & 0 & -1 & 0\\0 & 0 & 0 & -1\endVmatrix.
\mytag{5.25}
$$
For the components of the Levi-Civita connection in such 
a frame we have 
$$
\hskip -2em
\Gamma^k_{ij}=
\frac{c^{\,k}_{ij}}{2}
-\sum^3_{r=0}\sum^3_{s=0}\frac{c^{\,s}_{i\kern 0.5pt r}}{2}\,g^{kr}
\,g_{sj}-\sum^3_{r=0}\sum^3_{s=0}\frac{c^{\,s}_{j\kern 0.5ptr}}{2}
\,g^{kr}\,g_{s\kern 0.5pt i}.
\mytag{5.26}
$$
The formula \mythetag{5.26} is derived from the general formula
\thetag{6.3} in \mycite{9}. From the formulas \mythetag{5.25} and
\mythetag{5.26} we derive the following equalities:
$$
\xalignat 2
&\hskip -2em
L_{\boldsymbol\Upsilon_m}(g_{rj})=0,
&&\Gamma^k_{ij}-\Gamma^k_{j\kern 0.5pt i}=c^{\,k}_{ij}.
\mytag{5.27}
\endxalignat
$$
Moreover, from $\nabla_{\!m\,}g_{ij}=0$ in this case we derive
$$
\hskip -2em
\sum^3_{r=0}\sum^3_{s=0}g^{is}\,\Gamma^r_{ms}\,g_{rj}
=-\Gamma^i_{mj}.
\mytag{5.28}
$$
In the non-holonomic frame $\boldsymbol\Upsilon_0,\,\boldsymbol
\Upsilon_1,\,\boldsymbol\Upsilon_2,\,\boldsymbol\Upsilon_3$ the 
formulas \mythetag{5.19} are replaced by
$$
\hskip -2em
\aligned
&L_{\boldsymbol\Upsilon_s}(X^r)=\nabla_{\!s}X^r-
\sum^3_{m=0}\Gamma^r_{sm}\,X^m,\\
&L_{\boldsymbol\Upsilon_j}(X^i)=\nabla_{\!j}X^i-
\sum^3_{m=0}\Gamma^i_{jm}\,X^m.
\endaligned
\mytag{5.29}
$$
Applying \mythetag{5.27}, \mythetag{5.28}, and \mythetag{5.29}
to \mythetag{5.15}, we get the formula
$$
\hskip -2em
\gathered
\bundleboldX_K=\sum^3_{i=0}\sum^3_{m=0}X^m\,\varUpsilon^i_m
\,\frac{\partial}{\partial x^i}
+\sum^3_{i=0}\sum^3_{j=0}\!\left(\frac{1}{2}\nabla_{\!j}X^i
\,-\vphantom{\frac{1}{2}\sum^3_{r=0}
\sum^3_{s=0}g^{is}\,\nabla_{\!s}X^r\,g_{rj}}\right.\\
\left.-\,\frac{1}{2}\sum^3_{r=0}
\sum^3_{s=0}g^{is}\,\nabla_{\!s}X^r\,g_{rj}\,
-\sum^3_{m=0}X^m\,\Gamma^i_{mj}\!\right)\!v^j
\,\frac{\partial}{\partial v^i}.
\endgathered
\mytag{5.30}
$$
Similarly, applying \mythetag{5.29} to \mythetag{5.2} and taking 
into account \mythetag{5.27}, we get 
$$
\bundleboldX_N=\sum^3_{i=0}\sum^3_{m=0}X^m\,\varUpsilon^i_m
\,\frac{\partial}{\partial x^i}
+\sum^3_{i=0}\sum^3_{j=0}\!\left(\!\nabla_{\!j}X^i
-\sum^3_{m=0}X^m\,\Gamma^i_{mj}\!\right)\!v^j\,\frac{\partial}
{\partial v^i}.\quad
\mytag{5.31}
$$
The formulas \mythetag{5.30} to \mythetag{5.31} coincide with
\mythetag{5.21} to \mythetag{5.22} respectively, though 
$\Gamma^i_{mj}$ are not symmetric with respect to $m$ and $j$
in this case. \pagebreak The formulas \mythetag{5.30} and 
\mythetag{5.31} again lead to the formula \mythetag{5.23}, 
which is valid irrespective to the choice of a holonomic or 
non-holonomic frame in $M$.\par
     Note that the formulas \mythetag{5.29} are the same as the
formula \thetag{3.10} in \mycite{10}. The formula \mythetag{5.29}
resembles the formulas \thetag{3.11} and \thetag{3.12} in
\mycite{10}. However, it doesn't coincide with them. Unlike
\mycite{10} and \mycite{11}, for the sake of simplicity in 
this paper I do not use principal fiber bundles at all.\par
\head
6. Commutation relationships for Kosmann-Lie derivatives.
\endhead
     Regular Lie derivatives acting upon tensorial fields form
a representation of the Lie algebra of vector fields in $M$. This
fact is expressed by the formula
$$
\hskip -2em
[L_{\bold X},L_{\bold Y}]=L_{[\bold X,\bold Y]}.
\mytag{6.1}
$$
Commutation relationships for Kosmann-Lie derivatives are different
from \mythetag{6.1}: 
$$
\hskip -2em
[\Cal L_{\bold X},\Cal L_{\bold Y}]=\Cal L_{[\bold X,\bold Y]}
+S_{\bold X,\bold Y}.
\mytag{6.2}
$$
Here $S_{\bold X,\bold Y}$ is a degenerate differentiation given 
by the tensor field $\bold S_{\bold X,\bold Y}$, where 
$$
\hskip -2em
\bold S_{\bold X,\bold Y}=L_{\bold X}(\bold S_{\bold Y})
-L_{\bold Y}(\bold S_{\bold X})-\bold S_{[\bold X,\bold Y]}
+[\bold S_{\bold X},\bold S_{\bold Y}].
\mytag{6.3}
$$
By means of direct calculations the formula \mythetag{6.3} can be 
reduced to 
$$
\hskip -2em
\bold S_{\bold X,\bold Y}=-[\bold S_{\bold X},\bold S_{\bold Y}].
\mytag{6.4}
$$
Applying \mythetag{6.4} to \mythetag{6.2}, we get 
$$
\hskip -2em
[\Cal L_{\bold X},\Cal L_{\bold Y}]=\Cal L_{[\bold X,\bold Y]}
-[\bold S_{\bold X},\bold S_{\bold Y}].
\mytag{6.5}
$$
The commutator $[\bold S_{\bold X},\bold S_{\bold Y}]$ in the formulas 
\mythetag{6.3}, \mythetag{6.4}, and \mythetag{6.5} is understood as a 
commutator of two operator-valued tensorial fields:
$$
\hskip -2em
[\bold S_{\bold X},\bold S_{\bold Y}]=\bold S_{\bold X}\compos
\bold S_{\bold Y}-\bold S_{\bold X}\compos\bold S_{\bold Y}.
\mytag{6.6}
$$
In the coordinate representation the commutator \mythetag{6.6} turns
to the commutator of two matrices whose components are calculated
according to the formula \mythetag{5.24}.\par
     Note that in general case the commutator \mythetag{6.6} is not zero.
For this reason $[\Cal L_{\bold X},\Cal L_{\bold Y}]\neq\Cal L_{[\bold X,
\bold Y]}$. This fact is pointed out in \mycite{10}.
\head
7. Kosmann-Lie derivatives for Weyl spinors.
\endhead     
     The bundle of Weyl spinors is a two-dimensional complex vector bundle 
over the space-time manifold $M$. We denote it $SM$. The spinor bundle $SM$
is related to the tangent bundle $TM$ in some special way. The relation of
$SM$ and $TM$ is formulated in terms of frames. It is based on the well-known 
group homomorphism
$$
\hskip -2em
\phi:\ \MatGrSL(2,\Bbb C)\to\MatGrSO^+(1,3,\Bbb R).
\mytag{7.1}
$$
Let $\boldsymbol\Upsilon_0,\,\boldsymbol\Upsilon_1,\,\boldsymbol\Upsilon_2,
\,\boldsymbol\Upsilon_3$ be a positively polarized right orthonormal frame
of the tangent bundle $TM$. By definition (see Section~5 in \mycite{12} or
Section 1 in \mycite{13}) it is cano\-nically associated with some frame 
$\boldsymbol\Psi_1,\,\boldsymbol\Psi_2$ of $SM$ in such a way that if two 
positively polarized right orthonormal frames $\boldsymbol\Upsilon_0,\,
\boldsymbol\Upsilon_1,\,\boldsymbol\Upsilon_2,\,\boldsymbol\Upsilon_3$ and 
$\tilde{\boldsymbol\Upsilon}_0,\,\tilde{\boldsymbol\Upsilon}_1,
\,\tilde{\boldsymbol\Upsilon}_2,\,\tilde{\boldsymbol\Upsilon}_3$ are
bound with the transition matrices $S$ and $T=S^{-1}$ in the formulas 
$$
\xalignat 2
&\tilde{\boldsymbol\Upsilon}_{\!i}=\sum^3_{j=0}S^j_i\,\boldsymbol
\Upsilon_{\!j},
&&\boldsymbol\Upsilon_{\!i}=\sum^3_{j=0}T^j_i\,\tilde{\boldsymbol
\Upsilon}_{\!j},
\endxalignat
$$
then their associated spinor frames $\boldsymbol\Psi_1,\,\boldsymbol\Psi_2$ 
and $\tilde{\boldsymbol\Psi}_1,\,\tilde{\boldsymbol\Psi}_2$ are bound with 
the transition matrices $\goth S$ and $\goth T=\goth S^{-1}$ in the formulas 
$$
\xalignat 2
&\tilde{\boldsymbol\Psi}_{\!i}=\sum^2_{j=1}\goth S^j_i\,
\boldsymbol\Psi_{\!j},
&&\boldsymbol\Psi_{\!i}=\sum^2_{j=1}\goth T^j_i\,
\tilde{\boldsymbol\Psi}_{\!j}
\endxalignat
$$
and the spacial transition matrices $S$ and $T$ are produced from 
the spinor transition matrices $\goth S$ and $\goth T$ through the 
group homomorphism \mythetag{7.1}:
$$
\xalignat 2
&S=\phi(\goth S),
&&T=\phi(\goth T).
\endxalignat
$$\par
     A spinor frame $\boldsymbol\Psi_1,\,\boldsymbol\Psi_2$ canonically
associated with some positively polarized right orthonormal tangent frame
$\boldsymbol\Upsilon_0,\,\boldsymbol\Upsilon_1,\,\boldsymbol\Upsilon_2,
\,\boldsymbol\Upsilon_3$ is called an {\it orthonormal frame\/}
of the spinor bundle $SM$. In order to visualize this canonical frame 
association in $SM$ and $TM$ we use the following diagram:
$$
\hskip -2em
\aligned
&\boxit{\lower 5pt\hbox{Orthonormal frames}}{}\to
\boxit{Positively polarized}{right orthonormal frames}\ .
\endaligned\quad
\mytag{7.2}
$$
There are two basic spin-tensorial fields in the bundle of Weyl spinors:
$$
\vcenter{\hsize 10cm
\offinterlineskip\settabs\+\indent
\vtrule
\hskip 1.2cm &\vtrule 
\hskip 5.2cm &\vtrule 
\hskip 2.8cm &\vtrule 
\cr\hrule 
\+\vtrule
\hfill\,Symbol\hfill&\vtrule
\hfill Name\hfill &\vtrule
\hfill Spin-tensorial\hfill &\vtrule\cr
\vskip -0.2cm
\+\vtrule
\hfill &\vtrule
\hfill \hfill&\vtrule
\hfill type\hfill&\vtrule\cr\hrule
\+\vtrule
\hfill $\bold d$\hfill&\vtrule
\hfill Skew-symmetric metric tensor\hfill&\vtrule
\hfill $(0,2|0,0|0,0)$\hfill&\vtrule\cr\hrule
\+\vtrule
\hfill$\bold G$\hfill&\vtrule
\hfill Infeld-van der Waerden field\hfill&\vtrule
\hfill $(1,0|1,0|0,1)$\hfill&\vtrule\cr\hrule
}\ .\quad
\mytag{7.3}
$$
Their role for $SM$ is similar to that of the metric tensor $\bold g$
for $TM$. The spin-tensorial type in the table \mythetag{3.1} specifies
the number of indices in coordinate representation of fields. 
The first two numbers are the numbers of upper and lower spinor 
indices, the second two numbers are the numbers of upper and 
lower conjugate spinor indices, and the last two numbers are 
the numbers of upper and lower tensorial indices (they are also 
called spacial indices).\par
     The spin-tensorial fields $\bold d$ and $\bold G$ are introduced
by means of the explicit formulas for their components in canonically 
associated frame pairs \mythetag{7.2}. Let $\boldsymbol\Upsilon_0,
\,\boldsymbol\Upsilon_1,\,\boldsymbol\Upsilon_2$, $\boldsymbol\Upsilon_3$
be some positively polarized right orthonormal frame in $TM$ and let 
$\boldsymbol\Psi_1,\,\boldsymbol\Psi_2$ be its associated orthonormal
spinor frame in $SM$. The components of the Infeld-van der Waerden field
$\bold G$ \pagebreak in such a frame pair composed by two canonically 
associated frames are given by the following Pauli matrices:
$$
\xalignat 2
&\hskip -2em
G^{i\kern 0.5pt\bar i}_0=\Vmatrix 1 & 0\\ 0 & 1\endVmatrix
=\sigma_0,
&&G^{i\kern 0.5pt\bar i}_2=\Vmatrix 0 & -i\\ i & 0\endVmatrix
=\sigma_2,\\
\vspace{-1.4ex}
&&&\mytag{7.4}\\
\vspace{-1.4ex}
&\hskip -2em
G^{i\kern 0.5pt\bar i}_1=\Vmatrix 0 & 1\\ 1 & 0\endVmatrix
=\sigma_1,
&&G^{i\kern 0.5pt\bar i}_3=\Vmatrix 1 & 0\\ 0 & -1\endVmatrix
=\sigma_3.
\endxalignat
$$
The skew-symmetric metric tensor $\bold d$ is given by the matrix
$$
\hskip -2em
d_{ij}=\Vmatrix 0 & 1\\ 
\vspace{1ex} -1 & 0\endVmatrix
\mytag{7.5}
$$
in any orthonormal spinor frame $\boldsymbol\Psi_1,\,\boldsymbol\Psi_2$.
Unlike \mythetag{7.4}, the choice of the associated frame $\boldsymbol
\Upsilon_0,\,\boldsymbol\Upsilon_1,\,\boldsymbol\Upsilon_2,\,\boldsymbol
\Upsilon_3$ is inessential for the components of the matrix \mythetag{7.5}
since $\bold d$ has no spatial indices at all. The dual metric tensor
for $\bold d$ is given by the matrix 
$$
\hskip -2em
d^{ij}=\Vmatrix 0 & -1\\ 
\vspace{1ex} 1 & 0\endVmatrix.
\mytag{7.6}
$$
The matrix \mythetag{7.6} is inverse to the matrix \mythetag{7.5}.\par
     Let $\goth S$ be some matrix from the group $\MatGrSL(2,\Bbb C)$ and
let $S=\phi(\goth S)$ be its image under the homomorphism \mythetag{7.1}.
Then we have the relationships 
$$
\align
&\hskip -2em
\sum^2_{i=1}\sum^2_{\bar i=1}\goth S^a_i\,\sigma^{i\kern 0.5pt\bar i}_m
\,\overline{\goth S^{\raise 0.6pt\hbox{$\ssize\bar a$}}_{\bar i}}
=\sum^3_{k=1}S^k_m
\,\sigma^{a\kern 0.5pt\bar a}_{\,k},
\mytag{7.7}\\
&\hskip -2em
\sum^2_{i=1}\sum^2_{j=1}\goth S^i_a\,d_{ij}\,\goth S^j_b=d_{ab},
\mytag{7.8}
\endalign
$$
where $\sigma^{i\kern 0.5pt\bar i}_m$ and $\sigma^{a\kern 0.5pt\bar a}_{\,k}$ 
are the components of the Pauli matrices \mythetag{7.4}. In essential, the
relationships \mythetag{7.7} form a definition of the group homomorphism 
\mythetag{7.1} (see \mycite{13} for more details). The relationships 
\mythetag{7.8} are fulfilled due to $\goth S\in\MatGrSL(2,\Bbb C)$.\par
     Now let's take some arbitrary vector field $\bold X$ on $M$. Then 
$\bundleboldX_K$ is its Kosmann lifting to the tangent bundle $TM$ given by 
the formula \mythetag{5.30}. It induces local one-parametric group of local
diffeomorphisms
$$
\hskip -2em
\varphi_{\varepsilon}\!:\,TM\to TM
\mytag{7.9}
$$
that extends \mythetag{2.1}, but does not coincide with \mythetag{2.2}. The
Kosmann lifting $\bundleboldX_K$ of the vector field $\bold X$ is concordant
with the bundle structure of the tangent bundle $TM$ in the sense of the 
definition~\mythedefinition{4.1}. For this reason the local diffeomorphisms 
\mythetag{7.9} break into the series of linear mappings 
$$
\hskip -2em
\varphi_{\varepsilon}\!:\,T_p(M)\to T_q(M),
\mytag{7.10}
$$
where $q=\varphi_{\varepsilon}(p)$. According to the 
definition~\mythedefinition{5.1} the mappings \mythetag{7.10} are isometries.
Therefore, taking some positively polarized right orthonormal frame 
$\boldsymbol\Upsilon_0,\,\boldsymbol\Upsilon_1,\,\boldsymbol\Upsilon_2$,
$\boldsymbol\Upsilon_3$ in $TM$ and applying $\varphi_{\varepsilon}$ to
it, we get another orthonormal frame $\varphi_{\varepsilon}(\boldsymbol
\Upsilon_0)$, $\varphi_{\varepsilon}(\boldsymbol\Upsilon_1),
\,\varphi_{\varepsilon}(\boldsymbol\Upsilon_2),\,\varphi_{\varepsilon}
(\boldsymbol\Upsilon_3)$ in $TM$. Note that $\varphi_{\varepsilon}$ is
homotopic to the identical mapping. For this reason it preserves the 
discrete properties like polarization and orientation, i\.\,e\. 
$\varphi_{\varepsilon}(\boldsymbol\Upsilon_0),\,\varphi_{\varepsilon}
(\boldsymbol\Upsilon_1),\,\varphi_{\varepsilon}(\boldsymbol\Upsilon_2),
\,\varphi_{\varepsilon}(\boldsymbol\Upsilon_3)$ a positively polarized
right orthonormal frame. Assume that $\varepsilon$ is small enough so
that both points $p$ and $q=\varphi_{\varepsilon}(p)$ belong to the
domain of the frame $\boldsymbol\Upsilon_0,\,\boldsymbol\Upsilon_1,\,
\boldsymbol\Upsilon_2,\,\boldsymbol\Upsilon_3$. Then at the point $q$
we have
$$
\hskip -2em
\varphi_{\varepsilon}(\boldsymbol\Upsilon_{\!j})
=\sum^3_{i=0}V^i_j(\varepsilon)\,\boldsymbol\Upsilon_{\!i}.
\mytag{7.11}
$$
The Lorentzian matrix $V$ with the components $V^i_j(\varepsilon)\in
\MatGrSO^+(1,3,\Bbb R)$ in \mythetag{7.11} is a coordinate presentation 
of the linear mapping \mythetag{7.10} in the frame $\boldsymbol\Upsilon_0,
\,\boldsymbol\Upsilon_1,\,\boldsymbol\Upsilon_2,\,\boldsymbol\Upsilon_3$.
\par
    Each of the two positively polarized right orthonormal frames
$\boldsymbol\Upsilon_0,\,\boldsymbol\Upsilon_1,\,\boldsymbol\Upsilon_2,
\,\boldsymbol\Upsilon_3$ and $\varphi_{\varepsilon}(\boldsymbol
\Upsilon_0),\,\varphi_{\varepsilon}(\boldsymbol\Upsilon_1),
\,\varphi_{\varepsilon}(\boldsymbol\Upsilon_2),\,\varphi_{\varepsilon}
(\boldsymbol\Upsilon_3)$ is associated with some orthonormal frame
in $SM$. Using this fact, we define a linear mapping
$$
\hskip -2em
\varphi_{\varepsilon}\!:\,S_p(M)\to S_q(M)
\mytag{7.12}
$$
closing the following commutative diagram of frame associations:
$$
\hskip -2em
\CD
\boldsymbol\Psi_1,\,\boldsymbol\Psi_2 @>>>\boldsymbol\Upsilon_0,
\,\boldsymbol\Upsilon_1,\,\boldsymbol\Upsilon_2,\,\boldsymbol
\Upsilon_3\\
\red{@V\varphi_{\varepsilon}VV} @VV\varphi_{\varepsilon}V\\
\varphi_{\varepsilon}(\boldsymbol\Psi_1),\,\varphi_{\varepsilon}
(\boldsymbol\Psi_2)@>>>
\aligned
&\varphi_{\varepsilon}(\boldsymbol
\Upsilon_0),\,\varphi_{\varepsilon}(\boldsymbol\Upsilon_1),\\
\vspace{-2ex}
&\varphi_{\varepsilon}(\boldsymbol\Upsilon_2),\,\varphi_{\varepsilon}
(\boldsymbol\Upsilon_3)
\endaligned
\endCD
\mytag{7.13}
$$
In the coordinate form the linear mapping \mythetag{7.12}, which  closes
the diagram \mythetag{7.13},  is presented by some matrix $W\in
\MatGrSL(2,\Bbb C)$:
$$
\hskip -2em
\varphi_{\varepsilon}(\boldsymbol\Psi_{\!j})
=\sum^2_{i=1}W^i_{\!j}(\varepsilon)\,\boldsymbol\Psi_{\!i}.
\mytag{7.14}
$$
The horizontal arrows in the diagram \mythetag{7.13} are canonical
frame associations. For this reason the components of the matrices
$V$ and $W$ in \mythetag{7.11} and \mythetag{7.14} should satisfy
the relationships \mythetag{7.7} and \mythetag{7.8}:
$$
\align
&\hskip -2em
\sum^2_{i=1}\sum^2_{\bar i=1}W^a_i(\varepsilon)
\,\sigma^{i\kern 0.5pt\bar i}_m
\,\overline{W^{\raise 0.6pt\hbox{$\ssize\bar a$}}_{\bar i}
(\varepsilon)}
=\sum^3_{k=1}V^k_m(\varepsilon)
\,\sigma^{a\kern 0.5pt\bar a}_{\,k},
\mytag{7.15}\\
&\hskip -2em
\sum^2_{i=1}\sum^2_{j=1}W^i_a(\varepsilon)\ d_{ij}
\,W^j_b(\varepsilon)=d_{ab}.
\mytag{7.16}
\endalign
$$
According to the general recipe \mythetag{4.12}, now we pass from the 
matrices $W(\varepsilon)$ and $V(\varepsilon)$ to their expansions as 
$\varepsilon\to 0$, i\.\,e\. we write
$$
\hskip -2em
\aligned 
&W^i_{\!j}(\varepsilon)=\delta^i_j+W^i_{\!j}(x^0,x^1,x^2,x^3)\,\varepsilon
+\ldots,\\
\vspace{2ex}
&V^i_j(\varepsilon)=\delta^i_j+V^i_j(x^0,x^1,x^2,x^3)\,\varepsilon
+\ldots.
\endaligned 
\mytag{7.17}
$$
Note that the quantities $V^i_j$ in \mythetag{7.17} are already known.
They are taken from \mythetag{5.13}. In our special case, where 
$\boldsymbol\Upsilon_0,\,\boldsymbol\Upsilon_1,\,\boldsymbol\Upsilon_2,
\,\boldsymbol\Upsilon_3$ is an orthonormal frame, we can use a more simple 
formula for $V^i_j$. It is extracted from \mythetag{5.30}:
$$
\hskip -2em
V^i_j=-\frac{1}{2}\sum^3_{r=0}
\sum^3_{s=0}g^{is}\,\nabla_{\!s}X^r\,g_{rj}
+\frac{1}{2}\nabla_{\!j}X^i
-\sum^3_{m=0}X^m\,\Gamma^i_{mj}.
\mytag{7.18}
$$
As for the quantities $W^i_{\!j}$, they should be calculated by substituting 
\mythetag{7.17} back into \mythetag{7.15} and \mythetag{7.16}. Like in
the case of \mythetag{5.7}, we denote 
$$
\hskip -2em
W_{ij}=\sum^2_{s=1}W^s_i\,d_{sj}.
\mytag{7.19}
$$
Then, applying \mythetag{7.17} to \mythetag{7.16}, we derive the following 
formula:
$$
\hskip -2em
W_{ij}-W_{j\kern 0.4pt i}=0.
\mytag{7.20}
$$
Due to the formula \mythetag{7.20} the quantities \mythetag{7.19} are
symmetric with respect to the indices $i$ and $j$. Our next goal is to
resolve the relationships \mythetag{7.15} with respect to these
quantities $W_{ij}$.\par
     Let's substitute the expansions \mythetag{7.17} into \mythetag{7.15}
and collect the first order terms with respect to the parameter 
$\varepsilon$. As a result we get
$$
\hskip -2em
\sum^2_{i=1}W^a_i\,\sigma^{i\kern 0.5pt\bar a}_m
+\sum^2_{\bar i=1}\sigma^{a\kern 0.5pt\bar i}_m
\,\overline{W^{\raise 0.6pt\hbox{$\ssize\bar a$}}_{\bar i}}
=\sum^3_{k=0}V^k_m\,\sigma^{a\kern 0.5pt\bar a}_{\,k}.
\mytag{7.21}
$$
Keeping in mind that we use the frames $\boldsymbol\Upsilon_0,\,\boldsymbol
\Upsilon_1,\,\boldsymbol\Upsilon_2,\,\boldsymbol\Upsilon_3$ and $\boldsymbol
\Psi_1,\,\boldsymbol\Psi_2$ canonically associated to each other in the
sense of the diagram \mythetag{7.2}, we replace the components of Pauli 
matrices by the components of Infeld-van der Waerden field in \mythetag{7.21}:
$$
\hskip -2em
\sum^2_{i=1}W^a_i\,G^{i\kern 0.5pt\bar a}_m
+\sum^2_{\bar i=1}G^{a\kern 0.5pt\bar i}_m
\,\overline{W^{\raise 0.6pt\hbox{$\ssize\bar a$}}_{\bar i}}
=\sum^3_{k=0}V^k_m\,G^{a\kern 0.5pt\bar a}_k.
\mytag{7.22}
$$
In order to transform the formula \mythetag{7.22} we multiply it by 
$G^m_{u\bar u}$ and sum over the index $m$. Doing it, we use one of 
the identities
$$
\xalignat 2
&\hskip -2em
\sum^3_{m=0}G^{a\kern 0.5pt\bar a}_m\,G^{\kern 0.2pt m}_{u\kern 0.5pt\bar u}
=2\,\delta^a_u\,\delta^{\bar a}_{\bar u},
&&\sum^2_{u=1}\sum^2_{\bar u=1}G^{u\kern 0.5pt\bar u}_m
\,G^{\kern 0.2pt n}_{u\kern 0.5pt\bar u}=2\,\delta^n_m,
\qquad
\mytag{7.23}
\endxalignat
$$
where $G^{\kern 0.2pt m}_{u\kern 0.5pt\bar u}$ are the components of 
the inverse Infeld-van der Waerden field. They are produced from 
$G^{a\kern 0.5pt\bar a}_n$ by lowering upper spinor indices $a$ amd 
$\bar a$ and by raising lower spacial index $n$ according to the 
following formula:
$$
\hskip -2em
G^{\kern 0.2pt m}_{u\kern 0.5pt\bar u}=\sum^2_{a=1}\sum^2_{\bar a=1}
\sum^3_{n=0}G^{a\kern 0.5pt\bar a}_n\,d_{a\kern 0.3pt u}\,\bd_{\bar a
\kern 0.3pt\bar u}\,g^{nm}.
\mytag{7.24}
$$
The components of the conjugate spinor metric $\bbd$ in \mythetag{7.24} 
are produced from the components of $\bold d$ by means of the complex 
conjugation:
$$
\xalignat 2
&\bd_{\,\bar i\kern 0.5pt\bar j}=\overline{d_{\,\bar i\kern 0.5pt\bar j}},
&&\bd^{\kern 0.5pt\bar i\bar j}=\overline{d^{\raise 1.2pt \hbox{$\ssize\kern 
0.5pt\bar i\bar j$}}}.
\endxalignat
$$
The identities \mythetag{7.23} are taken from \mycite{14}. Applying them to 
\mythetag{7.22}, we get
$$
\hskip -2em
2\,W^a_u\,\delta^{\bar a}_{\bar u}+2\,\delta^a_u
\,\overline{W^{\raise 0.6pt\hbox{$\ssize\bar a$}}_{\bar u}}
=\sum^3_{k=0}\sum^3_{m=0}V^k_m\,G^{a\kern 0.3pt\bar a}_k
\,G^{\kern 0.2pt m}_{u\kern 0.5pt\bar u}.
\mytag{7.25}
$$
In order to use the symmetry \mythetag{7.20} we need to lower the indices 
$a$ and $\bar a$ in \mythetag{7.25}:
$$
\hskip -2em
W_{u\kern 0.3pt a}\,\bd_{\bar u\kern 0.3pt\bar a}+d_{ua}
\,\overline{W_{\bar u\kern 0.3pt \bar a}}
=\frac{1}{2}\sum^3_{k=0}\sum^3_{m=0}\sum^2_{s=1}\sum^2_{\bar s=1}
V^k_m\,G^{s\kern 0.3pt\bar s}_k\,d_{s\kern 0.3pt a}\,\bd_{\bar s
\kern 0.3pt\bar a}\,G^{\kern 0.2pt m}_{u\kern 0.5pt\bar u}.
\quad
\mytag{7.26}
$$
Note that $\overline{W_{\bar u\kern 0.3pt \bar a}}$ in \mythetag{7.26} 
is symmetric with respect to the indices $\bar u$ and $\bar a$, while 
$\bd^{\kern 0.5pt\bar a\kern 0.3pt\bar u}$ is skew-symmetric with respect 
to these indices. Therefore, we have
$$
\hskip -2em
\sum^2_{\bar u=1}\sum^2_{\bar a=1}\overline{W_{\bar u\kern 0.3pt \bar a}}
\ \bd^{\kern 0.5pt\bar a\kern 0.3pt\bar u}=0.
\mytag{7.27}
$$
Applying \mythetag{7.27}, we multiply \mythetag{7.26} by $\bd^{\kern 0.5pt
\bar a\kern 0.3pt\bar u}$ and sum up over the indices $\bar u$ and $\bar a$.
As a result we get the following formula for $W_{u\kern 0.3pt a}$:
$$
\hskip -2em
W_{u\kern 0.3pt a}=\frac{1}{4}\sum^3_{k=0}\sum^3_{m=0}
\sum^2_{s=1}\sum^2_{\bar s=1}G^{s\kern 0.3pt\bar s}_k\,V^k_m
\,G^{\kern 0.2pt m}_{u\kern 0.5pt\bar s}\,d_{s\kern 0.3pt a}.
\mytag{7.28}
$$\par 
     Now let's return back to the quantities $W^i_{\!j}$ by raising the index 
$a$ in \mythetag{7.28}. As a result of this standard procedure we obtain
$$
\hskip -2em
W^i_{\!j}=\frac{1}{4}\sum^3_{k=0}\sum^3_{m=0}
\sum^2_{\bar s=1}G^{i\kern 0.3pt\bar s}_k\,V^k_m
\,G^{\kern 0.2pt m}_{j\kern 0.3pt\bar s}.
\mytag{7.29}
$$
The next step is to substitute \mythetag{7.18} into \mythetag{7.29}. Doing 
it, let's recall that the metric connection $\Gamma$ has the unique extension 
$(\Gamma,\Alpha,\bar{\Alpha})$ to the spinor bundle $SM$. Its spinor 
components are given by the formula
$$
\hskip -2em
\gathered
\Alpha^i_{\kern 0.5pt r j}
=\sum^3_{k=0}\sum^3_{m=0}\sum^2_{\bar s=1}
\frac{G^{i\bar s}_k\,\Gamma^k_{rm}
\,G^{\,m}_{j\kern 0.5pt\bar s}}{4}\,-\\
-\sum^2_{\bar s=1}\sum^3_{q=0}
\frac{L_{\boldsymbol\Upsilon_{\!r}}\!(G^{i\bar s}_q)
\,G^{\,q}_{j\kern 0.5pt\bar s}}{4}
-\sum^2_{\bar i=1}\sum^2_{\bar j=1}
\frac{L_{\boldsymbol\Upsilon_{\!r}}\!(\bd_{\bar j\kern 0.5pt\bar i})
\,\bd^{\kern 0.5pt\bar i\bar j}\,\delta^{\,i}_j}{4}.
\endgathered
\mytag{7.30}
$$
The formula is derived in \mycite{14} and is verified in \mycite{15}. In 
our special case the frames $\boldsymbol\Upsilon_0,\,\boldsymbol\Upsilon_1,
\,\boldsymbol\Upsilon_2,\,\boldsymbol\Upsilon_3$ and $\boldsymbol\Psi_1,
\,\boldsymbol\Psi_2$ form a canonically associated pair in the sense of 
the diagram \mythetag{7.2}. In this case the formula \mythetag{7.30}
reduces to
$$
\hskip -2em
\Alpha^i_{\kern 0.5pt r j}
=\frac{1}{4}\sum^3_{k=0}\sum^3_{m=0}\sum^2_{\bar s=1}
G^{i\bar s}_{\!k}\,\Gamma^k_{rm}\,G^{\,m}_{\!j\kern 0.5pt\bar s}.
\mytag{7.31}
$$
The formulas \mythetag{7.31} and \mythetag{7.29} are very similar. Therefore,
substituting \mythetag{7.18} into \mythetag{7.29}, due to the presence of
$\Gamma^i_{mj}$ in \mythetag{7.18} we can write
$$
\hskip -2em
\gathered
W^i_{\!j}=-\frac{1}{8}\sum^3_{q=0}\sum^3_{p=0}\sum^3_{k=0}
\sum^3_{m=0}\sum^2_{\bar s=1}G^{i\bar s}_k\,g^{kp}\,\nabla_{\!p}X^q
\,g_{qm}\,G^{\kern 0.2pt m}_{j\kern 0.3pt\bar s}\,+\\
+\,\frac{1}{8}\sum^3_{k=0}\sum^3_{m=0}\sum^2_{\bar s=1}G^{i\bar s}_k
\,\nabla_{\!m}X^k\,G^{\kern 0.2pt m}_{j\kern 0.3pt\bar s}
-\sum^3_{m=0}X^m\Alpha^i_{mj}.
\endgathered
\mytag{7.32}
$$
Knowing the quantities \mythetag{7.18} and \mythetag{7.32} is sufficient to
construct a spinor extension of the Kosmann-Lie derivative $\Cal L_{\bold X}$.
Let $\bold Y$ be a spin-tensorial field of the type $(\varepsilon,\eta|\sigma,
\zeta|e,f)$. Then for the components of the field $\Cal L_{\bold X}(\bold Y)$
we have the formula
$$
\gathered
\Cal L_{\bold X}(\bold Y)^{a_1\ldots\,a_\varepsilon\bar a_1\ldots
\,\bar a_\sigma c_1\ldots\,c_e}_{b_1\ldots\,b_\eta\bar b_1\ldots
\,\bar b_\zeta d_1\ldots\,d_f}=
\sum^3_{m=0}X^m\,L_{\boldsymbol\Upsilon_m}(Y^{a_1\ldots\,a_\varepsilon\bar a_1
\ldots\,\bar a_\sigma c_1\ldots\,c_e}_{b_1\ldots\,b_\eta\ \bar b_1\ldots
\,\bar b_\zeta\ d_1\ldots\,d_f})\,-\\
-\sum^\varepsilon_{\mu=1}\sum^2_{v_\mu=1}W^{a_\mu}_{v_\mu}\ 
Y^{a_1\ldots\,v_\mu\,\ldots\,a_\varepsilon\bar a_1\ldots\,\bar a_\sigma
c_1\ldots\,c_e}_{\ b_1\ldots\,\ldots\,\ldots\,b_\eta\,\bar b_1\ldots\,
\bar b_\zeta\,d_1\ldots\,d_f}\,+\\
\kern 9em+\sum^\eta_{\mu=1}\sum^2_{w_\mu=1}W^{w_\mu}_{b_\mu}\
Y^{a_1\ldots\,\ldots\,\ldots\,a_\varepsilon\bar a_1\ldots\,\bar a_\sigma
c_1\ldots\,c_e}_{\,b_1\ldots\,w_\mu\,\ldots\,b_\eta\bar b_1\ldots\,
\bar b_\zeta d_1\ldots\,d_f}\,-\\
\kern -9em
-\sum^\sigma_{\mu=1}\sum^2_{v_\mu=1}
\overline{W^{\bar a_\mu}_{v_\mu}}\ 
Y^{a_1\ldots\,a_\varepsilon\bar a_1\ldots\,v_\mu\,\ldots\,\bar a_\sigma
c_1\ldots\,c_e}_{\,b_1\ldots\,b_\eta\,\bar b_1\ldots\,\ldots\,\ldots\,
\bar b_\zeta\,d_1\ldots\,d_f}\,+\\
\kern 9em +\sum^\zeta_{\mu=1}\sum^2_{w_\mu=1}
\overline{W^{w_\mu}_{\bar b_\mu}}\
Y^{a_1\ldots\,a_\varepsilon\bar a_1\ldots\,\ldots\,\ldots\,\bar a_\sigma
c_1\ldots\,c_e}_{\,b_1\ldots\,b_\eta\,\bar b_1\ldots\,w_\mu\,\ldots\,
\bar b_\zeta d_1\ldots\,d_f}\,-\\
\kern -9em-\sum^e_{\mu=1}\sum^3_{v_\mu=0}V^{c_\mu}_{v_\mu}\ 
Y^{a_1\ldots\,a_\varepsilon\bar a_1\ldots\,\bar a_\sigma
c_1\ldots\,v_\mu\,\ldots\,c_e}_{\,b_1\ldots\,b_\eta\,\bar b_1\ldots\,
\bar b_\zeta\,d_1\ldots\,\ldots\,\ldots\,d_f}\,+\\
\kern 9em+\sum^f_{\mu=1}\sum^3_{w_\mu=0}V^{w_\mu}_{b_\mu}\
Y^{a_1\ldots\,a_\varepsilon\bar a_1\ldots\,\bar a_\sigma
c_1\ldots\,\ldots\,\ldots\,c_e}_{\,b_1\ldots\,b_\eta\,\bar b_1\ldots
\,\bar b_\zeta\,d_1\ldots\,w_\mu\,\ldots\,d_f}.
\endgathered\qquad
\mytag{7.33}
$$
It is easy to see that \mythetag{7.33} is a version of \mythetag{4.14}. 
In the case of a purely tensorial field $\bold Y$, i\.\,e\. if 
$\varepsilon=0$, $\eta=0$, $\sigma=0$, and $\zeta=0$, the Kosmann-Lie 
derivative \mythetag{7.33} reduces to \mythetag{5.23}. However, 
in general case we cannot use this formula \mythetag{5.23} since the 
regular Lie derivative $L_{\bold X}$ has no spinor extension yet. For 
this reason, instead of the formula \mythetag{5.23}, in this case we 
write 
$$
\hskip -2em
\Cal L_{\bold X}=\nabla_{\bold X}+S_{\bold X}.
\mytag{7.34}
$$
Like in \mythetag{5.23}, by $S_{\bold X}$ in \mythetag{7.34} we denote a
degenerate differentiation. According to the results of \mycite{12}, each 
degenerate differentiation extended to spinors is defined by three 
spin-tensorial field of the types $(1,1|0,0|0,0)$, $(0,0|1,1|0,0)$, and 
$(0,0|0,0|1,1)$. We denote them $\eufb S_{\bold X}$, 
$\bar{\eufb S}_{\bold X}$, and $\bold S_{\bold X}$ respectively. Here are 
the components of $\bold S_{\bold X}$:
$$
\hskip -2em
S^{\kern 0.2pt i}_j(\bold X)=\frac{\nabla^iX_j-\nabla_{\!j}X^i}{2}.
\mytag{7.35}
$$
Comparing \mythetag{7.35} with \mythetag{5.24}, we see that 
$\bold S_{\bold X}$ in \mythetag{7.34} is different from that of 
\mythetag{5.23}. The formula \mythetag{7.35} is extracted from 
\mythetag{7.18}. Similarly, looking at \mythetag{7.32}, we find 
the components of the spin-tensorial field $\eufb S_{\bold X}$:
$$
\hskip -2em
\goth S^{\kern 0.2pt i}_j(\bold X)=
\sum^3_{k=0}\sum^3_{m=0}\sum^2_{\bar s=1}G^{i\bar s}_k
\,\frac{\nabla^k\!X_m-\nabla_{\!m}X^k}{8}
\,G^{\kern 0.2pt m}_{j\kern 0.3pt\bar s}.
\mytag{7.36}
$$
The components of $\bar{\eufb S}_{\bold X}$ are produced from \mythetag{7.36}
by means of the complex conjugation: $\bar{\goth S}^{\kern 0.2pt
\bar i}_{\bar j}(\bold X)=\overline{\goth S^{\kern 0.2pt
\bar i}_{\bar j}(\bold X)}$. For this components we derive the formula 
$$
\hskip -2em
\bar{\goth S}^{\kern 0.2pt\bar i}_{\bar j}(\bold X)=
\sum^3_{k=0}\sum^3_{m=0}\sum^2_{s=1}G^{s\bar i}_k
\,\frac{\nabla^k\!X_m-\nabla_{\!m}X^k}{8}\,G^{\kern 0.2pt m}_{s\kern 0.3pt
\bar j}.
\mytag{7.37}
$$
The formula \mythetag{7.34} complemented with \mythetag{7.35}, 
\mythetag{7.36}, and \mythetag{7.37} is equivalent to the formula 
\mythetag{7.33}.\par
      Let's consider a particular example of applying the formula 
\mythetag{7.33}. Assume that $\boldsymbol\psi$ is a spinor field, 
i\.\,e\. a field with the spin-tensorial type $(1,0|0,0|0,0)$. Then 
for the components of the spinor field $\Cal L_{\bold X}(\boldsymbol\psi)$ 
we have
$$
\Cal L_{\bold X}(\boldsymbol\psi)^i=\sum^3_{m=0}X^m\,\nabla_m\psi^i
+\sum^3_{k=0}\sum^3_{m=0}\sum^2_{s=1}\sum^2_{j=1}G^{s\bar i}_k
\,\frac{\nabla^k\!X_m-\nabla_{\!m}X^k}{8}\,G^{\kern 0.2pt m}_{s
\kern 0.3pt\bar j}\,\psi^j.
$$
This formula resembles the formula \thetag{3.19} in \mycite{10} and the 
formula \thetag{5.5${}^\prime$} in \mycite{11}.\par
     Note that the formulas \mythetag{7.18} and \mythetag{7.32} were 
derived under the assumption that $\boldsymbol\Upsilon_0,\,\boldsymbol
\Upsilon_1,\,\boldsymbol\Upsilon_2,\,\boldsymbol\Upsilon_3$ is a positively 
polarized right orthonormal frame in $TM$ and $\boldsymbol\Psi_1,
\,\boldsymbol\Psi_2$ its canonically associated orthonormal frame 
in $SM$. However, these formulas remain valid for an arbitrary frame 
pair provided we use the general formula 
$$
\gathered
\Gamma^k_{ij}=\sum^3_{r=0}\frac{g^{\kern 0.5pt kr}}{2}
\left(L_{\boldsymbol\Upsilon_{\!i}}\!(g_{rj})
+L_{\boldsymbol\Upsilon_{\!j}}\!(g_{i\kern 0.5pt r})
-L_{\boldsymbol\Upsilon_{\!r}}\!(g_{ij})\right)+\\
+\,\frac{c^{\,k}_{ij}}{2}
-\sum^3_{r=0}\sum^3_{s=0}\frac{c^{\,s}_{i\kern 0.5pt r}}{2}\,g^{kr}
\,g_{sj}-\sum^3_{r=0}\sum^3_{s=0}\frac{c^{\,s}_{j\kern 0.5ptr}}{2}
\,g^{kr}\,g_{s\kern 0.5pt i}.
\endgathered
$$
for $\Gamma^i_{mj}$ in \mythetag{7.18} instead of \mythetag{5.26} and the 
general formula \mythetag{7.30} for $\Alpha^i_{mj}$ in \mythetag{7.32} 
instead of \mythetag{7.31}. The formula \mythetag{7.33} is also valid for 
an arbitrary frame pair under the same provisions.
\head
8. Kosmann-Lie derivatives of the basic fields.
\endhead
     There are three basic field in the theory of Weyl spinors. Two of them 
$\bold d$ and $\bold G$ are listed in the table \mythetag{7.3}. The third is 
the metric tensor $\bold g$. Now we shall apply the Kosmann-Lie derivative
\mythetag{7.33} to these basic fields. For this purpose it is convenient
to choose some canonically associated pair of frames $\boldsymbol\Upsilon_0,
\,\boldsymbol\Upsilon_1,\,\boldsymbol\Upsilon_2,\,\boldsymbol\Upsilon_3$ and
$\boldsymbol\Psi_1,\,\boldsymbol\Psi_2$. In such a frame pair the components
of all basic fields are constants. Indeed, they are given by the formulas
\mythetag{5.25}, \mythetag{7.4}, and \mythetag{7.5}. Therefore
we have
$$
\align
&\sum^3_{m=0}X^m\,L_{\boldsymbol\Upsilon_m}(d_{ij})=0,
\mytag{8.1}\\
&\sum^3_{m=0}X^m\,L_{\boldsymbol\Upsilon_m}(G^{a\kern 0.2pt\bar a}_m)=0,
\mytag{8.2}\\
&\sum^3_{m=0}X^m\,L_{\boldsymbol\Upsilon_m}(g_{ij})=0.
\mytag{8.3}
\endalign
$$
Applying \mythetag{7.33} to $\bold d$ and taking into account 
\mythetag{8.1}, \mythetag{7.19}, and \mythetag{7.20}, we obtain 
$$
\hskip -2em
\Cal L_{\bold X}(\bold d)_{ij}=\sum^2_{s=1}W^s_i\,d_{sj}
+\sum^2_{s=1}W^s_j\,d_{i\kern 0.2pt s}=W_{ij}
-W_{j\kern 0.3pt i}=0.
\mytag{8.4}
$$
Similarly, applying \mythetag{7.33} to $\bold G$ and taking into 
account \mythetag{8.2} and \mythetag{7.22}, we get
$$
\Cal L_{\bold X}(\bold G)^{a\kern 0.2pt \bar a}_m=
-\sum^2_{i=1}W^a_i\,G^{i\kern 0.5pt\bar a}_m
-\sum^2_{\bar i=1}G^{a\kern 0.5pt\bar i}_m
\,\overline{W^{\raise 0.6pt\hbox{$\ssize\bar a$}}_{\bar i}}
+\sum^3_{k=0}V^k_m\,G^{a\kern 0.5pt\bar a}_k=0.
\mytag{8.5}
$$
And finally we apply the formula \mythetag{7.33} to the metric tensor
$\bold g$. As a result, taking into account \mythetag{8.3}, \mythetag{5.7},
and \mythetag{5.8}, we derive 
$$
\hskip -2em
\Cal L_{\bold X}(\bold g)_{ij}=\sum^3_{r=0}V^r_i\,g_{rj}
+\sum^3_{r=0}V^r_j\,g_{i\kern 0.4pt r}=V_{ij}+V_{j\kern 0.3pt i}=0.
\mytag{8.6}
$$
The formulas \mythetag{8.4}, \mythetag{8.5}, and \mythetag{8.6} are
summarized in the following theorem.
\mytheorem{8.1} For any vector field $\bold X$ in $M$ the basic tensorial
and spin-tensorial fields $\bold g$, $\bold d$, and $\bold G$ associated
with the bundle of Weyl spinors $SM$ are constant with respect to the
Kosmann-Lie derivative $\Cal L_{\bold X}$.
\endproclaim
\head 
9. Some concluding remarks.
\endhead 
     Note that the quantities $V^i_j$ for \mythetag{7.33} are taken from 
\mythetag{7.18}. However, they could be taken from \mythetag{5.1} either.
In the latter case the equality $\Cal L_{\bold X}(\bold d)=0$ would be 
preserved, but the equality $\Cal L_{\bold X}(\bold g)=0$ would be replaced 
by
$$
\Cal L_{\bold X}(\bold g)=L_{\bold X}(\bold g).
$$
As for the formula \mythetag{8.5}, it would be replaced by the following one:
$$
\Cal L_{\bold X}(\bold G)^{a\kern 0.2pt \bar a}_m=
\sum^3_{k=0}\frac{\nabla_{\!m}X^k+\nabla^k\!X_m}{2}
\ G^{a\kern 0.5pt\bar a}_k.
$$
This choice of $V^i_j$ is preferred in \mycite{11}. As for our choice of
$V^i_j$ in this paper, in \mycite{11} it is referred to as the ``metric Lie
derivative'' \pagebreak introduced by Bourguignon and Gauduchon in 
\mycite{16}. Since there are various approaches, I should regretfully 
conclude that there is no canonical definition of the Lie derivative for 
spinors thus far.

\enddocument
\end